\newif\iflongVersion
\longVersiontrue

\newif\ifoneColumn
\oneColumnfalse

\ifoneColumn
\documentclass[12pt,draftcls,onecolumn]{IEEEtran}
\else
\documentclass[journal]{IEEEtran}
\fi

\ifCLASSINFOpdf
\else
\fi
\IEEEoverridecommandlockouts                              
\usepackage{epsfig} 
\usepackage{amsmath} 
\usepackage{amssymb}  

\newtheorem{thm}{Theorem}
\newtheorem{prp}{Proposition}
\newtheorem{cor}{Corollary}

\newtheorem{alg}{Algorithm}
\newtheorem{lem}{Lemma}
\newtheorem{rem}{Remark}
\newtheorem{defin}{Definition}
\newenvironment{pf}{\smallbreak\noindent{\it Proof. }}{\hfill$\Box$\smallbreak}

{\begin{list}{}%
         {\setlength{\leftmargin}{#1}}%
         \item[]%
}{\end{list}}

\makeatletter
\renewcommand*\env@matrix[1][*\c@MaxMatrixCols c]{%
  \hskip -\arraycolsep
  \let\@ifnextchar\new@ifnextchar
  \array{#1}}
\makeatother

\begin{document}
%

\title{\LARGE \bf
On feasibility, stability and performance in \\
distributed model predictive control
}

%
%
%

\author{{\bf{Pontus Giselsson}} and {\bf{Anders Rantzer}}\\
Department of Automatic Control\\
Lund University\\
Box 118, SE-221 00 Lund, Sweden\\
{\small\texttt{\{pontusg,rantzer\}@control.lth.se}}}

\maketitle

\begin{abstract}

\iflongVersion
In distributed model predictive control (DMPC), where a centralized
optimization problem is solved in distributed fashion using dual
decomposition, it is important to keep the number of iterations in the
solution algorithm, i.e.
the amount of communication between subsystems, as small as possible.
At the same time, 
the number of iterations must be enough to give a feasible
solution to the optimization problem and to guarantee stability of the
closed loop system. In this paper, a stopping condition to the
distributed optimization algorithm that guarantees these
properties, is
presented. The stopping condition is based on two theoretical
contributions. First, since the optimization problem is solved using
dual decomposition, standard techniques to prove stability in model
predictive control (MPC), i.e. with a terminal
cost and a terminal constraint set that involve all state variables,
do not apply. For the case without a terminal cost or a
terminal constraint set, we present a new method to quantify the control
horizon needed to ensure stability and a prespecified performance.
Second, the stopping condition is based on a novel adaptive
constraint tightening approach. Using this adaptive constraint
tightening approach, we guarantee that a
primal feasible solution to the optimization problem is found and that
closed loop stability and performance is obtained. 
Numerical examples show that the number of iterations needed to
guarantee feasibility of the optimization problem, stability and a
prespecified performance of the closed-loop system can 
be reduced significantly using the proposed stopping condition.
\fi

\end{abstract}


\begin{IEEEkeywords}
Distributed model predictive control, performance guarantee,
stability, feasibility
\end{IEEEkeywords}

%


\section{Introduction}
\label{sec:introduction}
Distributed model predictive control (DMPC) can be divided into two main
categories. In the first category, local optimization problems that
are solved sequentially and that take
neighboring interaction and solutions into account, are solved in each
subsystem. This is done in
 \cite{richards} for
linear systems and in \cite{dunbar} for nonlinear
systems. In \cite{hermans} a DMPC scheme is presented in which stability is
proven by adding a constraint to the optimization problem that
requires a reduction of an explicit control Lyapunov function. In
\cite{jia,camponogara} stability is guaranteed
for systems satisfying a certain matching condition and if the
coupling interaction is small enough. In the second category, to which
the current paper belong, a centralized optimization problem with a 
sparse structure is solved using a distributed optimization algorithm.
This approach is taken in \cite{venkat2008} where stability is
guaranteed in every algorithm 
iteration. A drawback to this method is that full model
knowledge is assumed in each node. Other approaches in the DMPC
literature rely on dual decomposition to solve the centralized
\iflongVersion
MPC
\fi
problem in distributed fashion. This approach is taken in, e.g.
\cite{Negenborn07,wakasa,doan2010}, where a (sub)gradient algorithm
is used to solve 
the dual problem and in \cite{necoaraTAC} where the
algorithm is based on the 
smoothing technique presented in \cite{Nesterov2005}. 
Among these, the only
stability proof is given in \cite{doan2009,doan2010}, where a
terminal point constraint is set to the origin, which is very restrictive.

One reason for the lack of stability results in DMPC based on dual
decomposition, is that the standard techniques to prove stability in
MPC do not apply. In MPC, terminal costs and terminal constraint sets
that involve all
state variables are used to show stability of the closed loop system,
see \cite{mayne,rawlingsBook}. This is not compatible with dual
decomposition. However, results for stability in MPC without a
terminal constraint set or a terminal set, which fits also the DMPC
framework used here, are available \cite{grimm,grune}. In
\cite{grune}, a method to quantify the minimal control horizon that
guarantees stability and a prespecified performance is presented. This is based on
relaxed dynamic programming \cite{lincoln,gruneRantzer} and a
controllability assumption on the stage costs. In the current paper,
we take a similar approach to quantify the control horizon needed to
guarantee stability and a prespecified performance. The advantages of
our approach over the one in \cite{grune} are twofold; we can, by
solving a mixed integer linear program (MILP),
verify our controllability assumption, further we get an explicit
expression that relates the parameter in the controllability
assumption with the obtained closed loop performance.

Besides the stability result, the main contribution of this paper is a
stopping condition for DMPC controllers that use a 
distributed optimization algorithm based on dual decomposition. We
use the distributed algorithm presented in \cite{gisAutomatica}, but any
duality-based distributed algorithm, such as the standard dual ascent or ADMM
\cite{BoydDistributed}, can be used. These duality based algorithms
suffer from that primal feasibility is only guaranteed in the limit of
iterations. Constraint tightening, which was originally proposed for
robust MPC in \cite{Gossner},
can also be used to generate feasible solutions within finite number of
iterations, see \cite{DoanConstrTight}. However, the introduction of
constraint tightening complicates
stability analysis since the optimal value
function without constraint tightening is used to show stability, while
the optimization is performed with tightened constraints. 
This problem is addressed in \cite{DoanConstrTight} by
assuming that the difference between the optimal value
functions with and without constraint tightening is bounded by a
constant. However, to actually
compute such a constant is very difficult.
The stopping condition in this paper is based on a novel adaptive
constraint tightening approach that ensures feasibility w.r.t. the
original constraint set with a finite number of algorithm iterations.
In addition, the amount of constraint tightening is adapted until the 
difference between the
optimal value functions with and without constraint tightening
is bounded by a certain amount. This adaptation
makes it possible to guarantee, besides feasibility of the
optimization problem, also
stability of the closed-loop system, without stating additional,
unquantifiable assumptions.

\iflongVersion
The paper is organized as follows. In Section~\ref{sec:probForm} we
introduce the problem and present the distributed optimization
algorithm in \cite{gisAutomatica}.
In Section~\ref{sec:stopConditions} the 
stopping condition is presented and feasibility, stability, and
performance is analyzed. 
Section~\ref{sec:ctrbVer} is devoted to
computation of a controllability parameter in the controllability assumption.
A numerical example that shows
the efficiency of the proposed stopping condition, is presented in
Section~\ref{sec:numerical}. Finally, in Section \ref{sec:conclusions}
we conclude the paper.
\fi


\section{Problem setup and preliminaries}
\label{sec:probForm}
We consider linear dynamical systems of the form
\begin{align}
x_{t+1} &= Ax_{t}+Bu_{t},& x_{0} &= \bar{x}
\label{eq:plant}
\end{align}
where $x_t\in\mathbb{R}^n$ and $u_t\in\mathbb{R}^m$ denote the state and
control vectors at time $t$ and the pair $(A,B)$ is
assumed controllable. We introduce the
following state and control variable partitions
\begin{align}
\label{eq:xPartition} x_t&=[(x_t^1)^T,(x_t^2)^T,\ldots,(x_t^{M})^T]^T,\\
\label{eq:uPartition} u_t&=[(u_t^1)^T,(u_t^2)^T,\ldots,(u_t^{M})^T]^T
\end{align}
where the local variables $x_t^i\in\mathbb{R}^{n_i}$ and
$u_t^i\in\mathbb{R}^{m_i}$.
The $A$ and $B$ matrices are partitioned accordingly
\begin{align*}
A &= \begin{pmatrix}
A_{11} & \cdots & A_{1M}\\
\vdots & \ddots & \vdots \\
A_{M1} & \cdots & A_{MM}\\
\end{pmatrix}, & 
B &= \begin{pmatrix}
B_{11} & \cdots & B_{1M}\\
\vdots & \ddots & \vdots \\
B_{M1} & \cdots & B_{MM}\\
\end{pmatrix}.
\end{align*}
These matrices are assumed to have a sparse structure, i.e., some
$A_{ij}=0$ and $B_{ij}=0$ and the neighboring interaction is defined by
the following sets
\begin{equation*}
\mathcal{N}_i = \{j\in\{1,\ldots,M\}~|~{\hbox{if }} A_{ij} \neq 0 {\hbox{ or }}
B_{ij}\neq 0\}
\end{equation*}
for $i=1,\ldots,M$. This gives the following local dynamics
\begin{align*}
x^i_{t+1} &= \sum_{j\in\mathcal{N}_i} \left(A_{ij}x_t^j+B_{ij}u_t^{j}\right),& x_{0}^i &= \bar{x}_i
\end{align*}
for $i=1,\ldots,M$. The local control and state variables are
constrained, i.e., $u^i\in \mathcal{U}_i$ and $x^i\in \mathcal{X}_i$.
The constraint sets, $\mathcal{X}_i$, $\mathcal{U}_i$ are assumed to be
bounded polytopes containing zero in their respective interiors and
can hence be represented as
\ifoneColumn
\begin{align*}
\mathcal{X}_i &= \{x^i\in\mathbb{R}^{n_i}\:|\:C_x^ix^i \leq d_x^i \},&
\mathcal{U}_i &= \{u^i\in\mathbb{R}^{m_i}\:|\:C_u^iu^i \leq d_u^i \}
\end{align*}
\else
\begin{align*}
\mathcal{X}_i &= \{x^i\in\mathbb{R}^{n_i}\:|\:C_x^ix^i \leq d_x^i \},\\
\mathcal{U}_i &= \{u^i\in\mathbb{R}^{m_i}\:|\:C_u^iu^i \leq d_u^i \}
\end{align*}
\fi
where $C_x^i\in\mathbb{R}^{n_{c_{x^i}}\times n_i}$,
$C_u^i\in\mathbb{R}^{n_{c_{u^i}}\times m_i}$,
$d_x^i\in\mathbb{R}_{>0}^{n_{c_{x^i}}}$ and
$d_u^i\in\mathbb{R}_{>0}^{n_{c_{u^i}}}$. We also denote the total number of
linear inequalities describing all constraint sets by $n_{c} :=
\sum_{i=1}^M \left(n_{c_{x^i}}+n_{c_{u^i}}\right)$.
The global constraint sets are defined from the local ones through
\begin{align*}
\mathcal{X}&=\mathcal{X}_1\times\ldots\times\mathcal{X}_M, & \mathcal{U}=\mathcal{U}_1\times\ldots\times\mathcal{U}_M.
\end{align*}
We use a separable quadratic stage cost
\begin{equation*}
\ell(x,u) = \sum_{i=1}^M \ell_i(x^i,u^i) =
\frac{1}{2}\left(\sum_{i=1}^M (x^i)^TQ_ix^i+(u^i)^TR_iu^i\right)
\end{equation*}
where $Q_i\in\mathbb{S}_{++}^{n_i}$ and
$R_i\in\mathbb{S}_{++}^{m_i}$ for $i=1,\ldots,M$ and
$\mathbb{S}_{++}^{n}$ denotes the set of symmetric
positive definite matrices in $\mathbb{R}^{n\times n}$. The
optimal infinite horizon
cost from initial state $\bar{x}\in\mathcal{X}$ is defined by
\begin{equation}
V_{\infty}(\bar{x}) := \begin{tabular}[t]{rl}
$\displaystyle\min_{x,u}$ & $\displaystyle\sum_{t=0}^{\infty}\ell(x_{t},u_{t})$\\
s.t. & $x_{t}\in \mathcal{X}$\phantom{aa},\phantom{aa}$u_{t} \in \mathcal{U}$\\
& $x_{t+1} = Ax_{t}+Bu_{t}$\\
& $x_{0} = \bar{x}$.
\end{tabular}
\label{eq:infValue}
\end{equation}
Such infinite horizon optimization problems are in general intractable to
solve exactly. A common approach is to solve the problem approximately in receding horizon
fashion. To this end, we introduce the predicted state and control
sequences $\{z_{\tau}\}_{\tau=0}^{N-1}$ and
$\{v_{\tau}\}_{\tau=0}^{N-1}$ and the corresponding
stacked vectors
\begin{align*}
\mathbf{z}&=[z_{0}^T,\ldots,z_{N-1}^T]^T,&
\mathbf{v}&=[v_{0}^T,\ldots,v_{N-1}^T]^T
\end{align*}
where $z_{\tau}$ and $v_{\tau}$ are predicted states and controls
$\tau$ time steps ahead.
The predicted state and control variables
$z_{\tau}$, $v_{\tau}$ are partitioned into local
variables as in
\eqref{eq:xPartition} and \eqref{eq:uPartition} respectively. We also
introduce the following stacked local vectors
\begin{align*}
\mathbf{z}_i&=[(z_{0}^i)^T,\ldots,(z^i_{N-1})^T]^T,&
 \mathbf{v}_i&=[(v_{0}^i)^T,\ldots,(v^i_{N-1})^T]^T.
\end{align*}
Further, we introduce the tightened state and control constraint sets
\begin{align*}
(1-\delta)\mathcal{X}_i &= \{x^i\in\mathbb{R}^{n_i}\:|\:C_x^ix^i \leq (1-\delta)d_x^i \},\\
(1-\delta)\mathcal{U}_i &= \{u^i\in\mathbb{R}^{m_i}\:|\:C_u^iu^i
\leq (1-\delta)d_u^i \}
\end{align*}
where $\delta\in(0,1)$ decides the amount of relative constraint tightening.
The following
optimization problem, which has neither a
terminal cost nor a terminal constraint set, is solved in the DMPC
controller for the current state $\bar{x}\in\mathbb{R}^n$
\begin{equation}
V_{N}^{\delta}(\bar{x}) := \begin{tabular}[t]{rl}
$\displaystyle\min_{\mathbf{z}_{t},\mathbf{v}_{t}}$ & $\displaystyle\sum_{\tau=0}^{N-1}\ell(z_\tau,v_\tau)$\\
s.t. & $z_{\tau}\in (1-\delta)\mathcal{X}$, $\tau=0,\ldots,N-1$\\
& $v_{\tau} \in (1-\delta)\mathcal{U}$, $\tau=0,\ldots,N-1$\\
& $z_{\tau+1} = Az_{\tau}+Bv_{\tau}$, $\tau = 0,\ldots,N-2$\\
& $z_{0} = \bar{x}$.
\end{tabular}
\label{eq:truncValue}
\end{equation}
By stacking all decision variables into
one vector
\begin{equation}
\mathbf{y}=[z_0^T,\ldots,z_{N-1}^T,v_0^T,\ldots,v_{N-1}^T]^T\in\mathbb{R}^{(n+m)N}
\label{eq:chiDef}
\end{equation}
the optimization
problem \eqref{eq:truncValue} can more compactly be written as
\begin{equation}
V_{N}^{\delta}(\bar{x}) := \begin{tabular}[t]{rl}
$\displaystyle\min_{\mathbf{y}}$ &
$\frac{1}{2}\mathbf{y}^T \mathbf{H} \mathbf{y}$\\
s.t.  & $\mathbf{A}\mathbf{y}=\mathbf{b}\bar{x}$\\
& $\mathbf{C}\mathbf{y}\leq (1-\delta)\mathbf{d}$
\end{tabular}
\label{eq:truncValueCompact}
\end{equation}
where 
$\mathbf{H}\in\mathbb{S}_{++}^{(n+m)N},
\mathbf{A}\in\mathbb{R}^{n(N-1)\times(n+m)N},
\mathbf{b}\in\mathbb{R}^{n(N-1)\times n},
\mathbf{C}\in\mathbb{R}^{n_cN\times(n+m)N}$ and
$\mathbf{d}\in\mathbb{R}_{>0}^{Nn_c}$ are built accordingly.
Such sparse optimization problems can be solved in distributed
fashion using, e.g., the classical dual ascent, the alternating
direction of multipliers method (ADMM)
\cite{BoydDistributed}, or the recently developed algorithm in
\cite{gisAutomatica}. The algorithm in \cite{gisAutomatica} is a dual
accelerated gradient algorithm
and is used in the current paper for simplicity. Distribution of
these methods are enabled by
solving the dual problem to \eqref{eq:truncValue}.

The dual
problem to \eqref{eq:truncValueCompact} is created by introducing dual variables
$\boldsymbol{\lambda}\in\mathbb{R}^{n(N-1)}$ for the equality
constraints and dual variables $\boldsymbol{\mu}\in\mathbb{R}_{\geq
  0}^{Nn_c}$ for the 
inequality constraints. As shown in \cite{gisAutomatica}, the dual
problem can explicitly be written as
\ifoneColumn
\begin{align}
 \max_{\boldsymbol{\lambda},\boldsymbol{\mu}\geq 0} -\frac{1}{2}(\mathbf{A}^T\boldsymbol{\lambda}+\mathbf{C}^T\boldsymbol{\mu})^T\mathbf{H}^{-1}(\mathbf{A}^T\boldsymbol{\lambda}+\mathbf{C}^T\boldsymbol{\mu})-
\boldsymbol{\lambda}^T\mathbf{b}\bar{x}-\boldsymbol{\mu}^T\mathbf{d}(1-\delta)
\label{eq:dualProb}
\end{align}
\else
\begin{multline}
\label{eq:dualProb}\max_{\boldsymbol{\lambda},\boldsymbol{\mu}\geq 0} -\frac{1}{2}(\mathbf{A}^T\boldsymbol{\lambda}+\mathbf{C}^T\boldsymbol{\mu})^T\mathbf{H}^{-1}(\mathbf{A}^T\boldsymbol{\lambda}+\mathbf{C}^T\boldsymbol{\mu})-\\
-\boldsymbol{\lambda}^T\mathbf{b}\bar{x}-\boldsymbol{\mu}^T\mathbf{d}(1-\delta)
\end{multline}
\fi
and we define the minimand in \eqref{eq:dualProb} as the dual
function for initial condition
$\bar{x}\in\mathbb{R}^n$, i.e.,
\ifoneColumn
\begin{align}
 D_{N}^{\delta}(\bar{x},\boldsymbol{\lambda},\boldsymbol{\mu})
:=
-\frac{1}{2}(\mathbf{A}^T\boldsymbol{\lambda}+\mathbf{C}^T\boldsymbol{\mu})^T\mathbf{H}^{-1}(\mathbf{A}^T\boldsymbol{\lambda}+\mathbf{C}^T\boldsymbol{\mu})-
\boldsymbol{\lambda}^T\mathbf{b}\bar{x}-\boldsymbol{\mu}^T\mathbf{d}(1-\delta).
\label{eq:dualFcn}
\end{align}
\else
\begin{multline}
\label{eq:dualFcn} D_{N}^{\delta}(\bar{x},\boldsymbol{\lambda},\boldsymbol{\mu})
:=
-\frac{1}{2}(\mathbf{A}^T\boldsymbol{\lambda}+\mathbf{C}^T\boldsymbol{\mu})^T\mathbf{H}^{-1}(\mathbf{A}^T\boldsymbol{\lambda}+\mathbf{C}^T\boldsymbol{\mu})-\\
-\boldsymbol{\lambda}^T\mathbf{b}\bar{x}-\boldsymbol{\mu}^T\mathbf{d}(1-\delta).
\end{multline}
\fi
The distributed algorithm presented in \cite{gisAutomatica} that solves
\eqref{eq:truncValueCompact}, is a dual accelerated gradient method
described by the following global iterations
\ifoneColumn
\begin{align}
\label{eq:accGrad1}\mathbf{y}^k &= -\mathbf{H}^{-1}(\mathbf{A}^T\boldsymbol{\lambda}^k+\mathbf{C}^T\boldsymbol{\mu}^k)\\
\label{eq:accGrad2}\mathbf{\bar{y}}^k &= \mathbf{y}^k+\frac{k-1}{k+2}(
\mathbf{y}^k-\mathbf{y}^{k-1})\\
\label{eq:accGrad3}\boldsymbol{\lambda}^{k+1} &=
\boldsymbol{\lambda}^k+\frac{k-1}{k+2}(\boldsymbol{\lambda}^k-\boldsymbol{\lambda}^{k-1})+\frac{1}{L}(\mathbf{A}\mathbf{\bar{y}}^k-\mathbf{b}\bar{x})\\
\label{eq:accGrad4}
\boldsymbol{\mu}^{k+1} &=
\max\bigg(0,\boldsymbol{\mu}^k+\frac{k-1}{k+2}(\boldsymbol{\mu}^k-\boldsymbol{\mu}^{k-1})+
\frac{1}{L}(\mathbf{C}\mathbf{\bar{y}}^k-\mathbf{d}\left(1-\delta)\right)\bigg)
\end{align}
\else
\begin{align}
\label{eq:accGrad1}\mathbf{y}^k &= -\mathbf{H}^{-1}(\mathbf{A}^T\boldsymbol{\lambda}^k+\mathbf{C}^T\boldsymbol{\mu}^k)\\
\label{eq:accGrad2}\mathbf{\bar{y}}^k &= \mathbf{y}^k+\frac{k-1}{k+2}(
\mathbf{y}^k-\mathbf{y}^{k-1})\\
\label{eq:accGrad3}\boldsymbol{\lambda}^{k+1} &=
\boldsymbol{\lambda}^k+\frac{k-1}{k+2}(\boldsymbol{\lambda}^k-\boldsymbol{\lambda}^{k-1})+\frac{1}{L}(\mathbf{A}\mathbf{\bar{y}}^k-\mathbf{b}\bar{x})\\
\nonumber
\boldsymbol{\mu}^{k+1} &=
\max\bigg(0,\boldsymbol{\mu}^k+\frac{k-1}{k+2}(\boldsymbol{\mu}^k-\boldsymbol{\mu}^{k-1})+\\
\label{eq:accGrad4}&\qquad\qquad\qquad\qquad\quad+\frac{1}{L}(\mathbf{C}\mathbf{\bar{y}}^k-\mathbf{d}\left(1-\delta)\right)\bigg)
\end{align}
\fi
where $k$ is the iteration number and
$L=\|[\mathbf{A}^T,\mathbf{C}^T]^T\mathbf{H}^{-1}[\mathbf{A}^T,\mathbf{C}^T]\|$,
which is the Lipschitz constant to the gradient of the dual function
\eqref{eq:dualFcn}. The reader is referred to \cite{gisAutomatica} for
details on how to distribute the algorithm \eqref{eq:accGrad1}-\eqref{eq:accGrad4}.

\subsection{Notation}
We define 
$\mathbb{N}_{\geq T}$ the 
set of natural numbers $t \geq T$. The norm
$\|\cdot\|$ refers to the Euclidean norm or the induced Euclidean
norm unless otherwise is
specified and $\langle\cdot,\cdot\rangle$ refers to the inner product
in Euclidean space. The norm $\|x\|_M = \sqrt{x^TMx}$. 
The optimal state and control sequences to
\eqref{eq:truncValue} for initial value $x$ and constraint tightening
$\delta$ are denoted
$\{z^*_{\tau}(x,\delta)\}_{\tau=0}^{N-1}$ and
$\{v^*_{\tau}(x,\delta)\}_{\tau=0}^{N-1}$ respectively and the optimal solution
to the equivalent problem \eqref{eq:truncValueCompact} by
$\mathbf{y}^*(x,\delta)$. The state and control sequences for
iteration $k$ in \eqref{eq:accGrad1}-\eqref{eq:accGrad4} are
denoted $\{z^k_{\tau}(x,\delta)\}_{\tau=0}^{N-1}$ and
$\{v^k_{\tau}(x,\delta)\}_{\tau=0}^{N-1}$ respectively.
The initial state and constraint tightening
arguments $(x,\delta)$ are dropped when no ambiguities can arise. 

\subsection{Definitions and assumptions}

We
adopt the convention that $V_{N}^{\delta}(\bar{x})=\infty$ for states
$\bar{x}\in\mathbb{R}^n$ that result in \eqref{eq:truncValueCompact}
being infeasible.
We define by $\mathbb{X}_\infty$ the
set for which \eqref{eq:infValue} is feasible and
we define the minimum of the stage-cost
$\ell$ for fixed $x$ as
\begin{equation*}
\ell^*(x) := \min_{u\in\mathcal{U}} \ell(x,u) = \frac{1}{2}x^TQx.
\end{equation*}
Further, $\kappa$ is the smallest scalar such that 
$\kappa Q-A^TQA\succeq 0$.
The state sequence resulting from applying
$\{v_{\tau}\}_{\tau=0}^{N-1}$ to \eqref{eq:plant} is denoted by
$\{\xi_{\tau}\}_{\tau=0}^{N-1}$, i.e.,
\begin{align}
\label{eq:xiDef}\xi_{\tau+1} &= A\xi_{\tau}+Bv_{\tau}, & \xi_{0} &= \bar{x}.
\end{align}
We introduce $\boldsymbol{\xi} = [(\xi_{0})^T,\ldots,(\xi_{N-1})^T]^T$
and define the primal cost
\ifoneColumn
\begin{align}
\label{eq:primCost}
&  P_N(\bar{x},\mathbf{v}) 
:=\left\{\begin{array}{ll}
\displaystyle\sum_{\tau=0}^{N-1}\ell(\xi_\tau,v_\tau)&
{\hbox{if }} \boldsymbol{\xi}\in \mathcal{X}^N {\hbox{ and }}
 \mathbf{v}\in \mathcal{U}^N {\hbox{ and}}
{\hbox{ \eqref{eq:xiDef} holds}}\\
\infty & {\hbox{else}}
\end{array}
\right.
\end{align}
\else
\begin{align}
\label{eq:primCost}
&  P_N(\bar{x},\mathbf{v}) 
:=\left\{\begin{array}{ll}
\displaystyle\sum_{\tau=0}^{N-1}\ell(\xi_\tau,v_\tau)&
{\hbox{if }} \boldsymbol{\xi}\in \mathcal{X}^N {\hbox{ and }}
 \mathbf{v}\in \mathcal{U}^N\\& {\hbox{~~~and}}
{\hbox{ \eqref{eq:xiDef} holds}}\\
\infty & {\hbox{else}}
\end{array}
\right.
\end{align}
\fi
where $\mathcal{X}^N$ and $\mathcal{U}^N$ are the state and control
constraints for the full horizon. 
We also introduce the shifted control sequence $\mathbf{v}_s =
[(v_1)^T,\ldots,(v_{N-1})^T,0^T]^T$. We have
$P_N(\bar{x},\mathbf{v}^k) \geq V_N(\bar{x})$ and
$P_N(A\bar{x}+Bv_0^k,\mathbf{v}_s^k)\geq V_N(A\bar{x}+Bv_0^k)$ for every
algorithm iteration $k$. We denote by $\{\xi^k_{\tau}\}_{\tau=0}^{N-1}$ the state
sequence that satisfies \eqref{eq:xiDef} using controls
$\{v^k_{\tau}\}_{\tau=0}^{N-1}$. The definition of the cost
\eqref{eq:primCost} implies
\begin{equation}
\label{eq:primCostRelation} P_N(\bar{x},\mathbf{v}^k) = P_N(A\bar{x}+Bv_0^k,\mathbf{v}_s^k)+\ell(\bar{x},v_0^k)-\ell^*(A\xi_{N-1}^k)
\end{equation}
if $v_0^k\in\mathcal{U}$, $\bar{x}\in\mathcal{X}$ and $A\xi_{N-1}^k\in\mathcal{X}$.



\section{Stopping condition}
\label{sec:stopConditions}
Rather than finding the optimal solution in each time step in the
MPC controller, the most 
important task is to find a control action that
gives desirable closed loop properties such as stability, feasibility,
and a
desired performance. Such properties can sometimes be ensured well before
convergence to the optimal solution. To benefit from this observation, a stopping
condition is developed that allows the iterations to stop when the
desired performance, stability, and feasibility can be guaranteed.
Before the stopping condition is introduced,
we briefly go through the main ideas below.

\subsection{Main ideas}

The
distributed nature of the optimization algorithm makes it unsuitable
for centralized terminal costs and terminal constraints. Thus,
stability and performance need to be ensured without these
constructions. We define the following infinite horizon
performance for feedback control law $\nu$
\begin{equation}
\label{eq:perfMetric}V_{\infty,\nu}(\bar{x}) = \sum_{t=0}^\infty\ell(x_t,\nu(x_t))
\end{equation}
where $x_{t+1} = Ax_{t}+B\nu(x_{t})$ and $x_0=\bar{x}$. For a
given performance 
parameter $\alpha\in(0,1]$ and control law $\nu$, it is known (cf.
\cite{lincoln,gruneRantzer,grune}) that the following decrease in the 
optimal value function
\begin{equation}
\label{eq:relDyn} V_N^0(x_t) \geq V_N^0(Ax_t+B\nu(x_t))+\alpha\ell(x_t,\nu(x_t))
\end{equation}
for every $t\in\mathbb{N}_{\geq 0}$
gives stability and closed loop performance according to
\begin{equation}
\label{eq:perfRes}\alpha V_{\infty,\nu}(\bar{x}) \leq V_{\infty}(\bar{x}).
\end{equation}
Analysis of the control horizon $N$ needed for an MPC control law without terminal
cost and terminal constraints such that \eqref{eq:relDyn}
holds, 
is performed in \cite{gruneRantzer,grune} and also in this paper.
Once a control horizon $N$ is known such that \eqref{eq:relDyn} is
guaranteed, the performance
result \eqref{eq:perfRes} relies on computation of the optimal
solution to the MPC optimization problem in every time step. An exact
optimal solution cannot be computed and the idea behind this paper is
to develop stopping conditions that enable early termination of the
optimization algorithm with
maintained feasibility, stability, and performance guarantees.
The idea behind our stopping
condition is to compute a lower bound to
$V_N^0(x)$ through the dual function
$D_N^0(x,\boldsymbol{\lambda}^k,\boldsymbol{\mu}^k)$ and an upper
bound to the next step value function
$V_N^0(Ax+Bv_0^k)$ through a feasible solution 
$P_N(Ax+Bv_0^k,\mathbf{v}_s^k)$. If at iteration $k$ the following
test is satisfied 
\begin{equation}
\label{eq:perfTest} D_N^0(\bar{x},\boldsymbol{\lambda}^k,\boldsymbol{\mu}^k) \geq
P_N(A\bar{x}+Bv_0^k,\mathbf{v}_s^k)+\alpha\ell(\bar{x},v_0^k)
\end{equation}
the performance condition \eqref{eq:relDyn} holds since
\begin{align*}
V_N^0(\bar{x})&\geq D_N^0(\bar{x},\boldsymbol{\lambda}^k,\boldsymbol{\mu}^k) \geq
P_N(A\bar{x}+Bv_0^k,\mathbf{v}_s^k)+\alpha\ell(\bar{x},v_0^k)\\
&\geq V_N^0(A\bar{x}+Bv_0^k)+\alpha\ell(\bar{x},v_0^k).
\end{align*}
This implies that stability and the performance result
\eqref{eq:perfRes} can be guaranteed with finite algorithm iterations
$k$ by using control action $v_0^k$.

The test \eqref{eq:perfTest} includes computation of 
$P_N(A\bar{x}+Bv_0^k,\mathbf{v}_s^k)$ which is a feasible solution to the
optimization problem in the following step. A feasible solution cannot be expected with finite
number of iterations $k$ for duality-based methods since 
primal feasibility is only guaranteed in the limit of iterations.
Therefore we introduce 
tightened state and control constraint sets
$(1-\delta)\mathcal{X}$, $(1-\delta)\mathcal{U}$ with $\delta\in(0,1)$ and use these in the
optimization problem. By generating a
state trajectory $\{\xi_{\tau}^k\}_{\tau=0}^{N-1}$ from the
control trajectory $\{v_{\tau}^k\}_{\tau=0}^{N-1}$ that satisfies
the equality constraints \eqref{eq:xiDef}, we will see that
$\{\xi_{\tau}^k\}_{\tau=0}^{N-1}$ satisfies the original inequality
constraints with finite number of iterations. Thus, a primal feasible
solution $P_N(A\bar{x}+Bv_0^k,\mathbf{v}_s^k)$ can be generated after a finite
number of algorithm iterations $k$. However, since the
optimization now is performed over a tightened constraint set, the dual function
value $D_N^\delta(\bar{x},\boldsymbol{\lambda},\boldsymbol{\mu})$ is not a
lower bound to $V_N^0(\bar{x})$ and cannot be used directly in the test
\eqref{eq:perfTest} to
ensure stability and the performance specified by \eqref{eq:perfRes}.
In the following lemma we show a
relation between the dual function value when using 
the tightened constraint sets and the optimal value function when
using the original constraint sets.
\begin{lem}
For every $\bar{x}\in\mathbb{R}^n$, $\boldsymbol{\lambda}\in\mathbb{R}^{n(N-1)}$ and
$\boldsymbol{\mu}\in\mathbb{R}_{\geq 0}^{Nn_c}$ we have that
\begin{equation*}
V_N^0(\bar{x})\geq
D_N^{\delta}(\bar{x},\boldsymbol{\lambda},\boldsymbol{\mu})-\delta\boldsymbol{\mu}^T\mathbf{d}.
\end{equation*}
\label{lem:tighteningRelationsVal}
\end{lem}
\iflongVersion
\begin{pf}
From the definition of the dual function \eqref{eq:dualFcn} we get that
\begin{align*}
D_{N}^{\delta}(\bar{x},\boldsymbol{\lambda},\boldsymbol{\mu})&=
D_{N}^{0}(\bar{x},\boldsymbol{\lambda},\boldsymbol{\mu})+\delta\mathbf{d}^T\boldsymbol{\mu}.
\end{align*}
By weak duality we get
\begin{align*}
V_{N}^{0}(\bar{x})&\geq
D_{N}^{0}(\bar{x},\boldsymbol{\lambda},\boldsymbol{\mu})=
D_{N}^{\delta}(\bar{x},\boldsymbol{\lambda},\boldsymbol{\mu})-\delta\mathbf{d}^T\boldsymbol{\mu}.
\end{align*}
This completes the proof.
\end{pf}
\fi

The presented lemma enables computation of a lower bound to $V_N^0(\bar{x})$
at algorithm iteration $k$ that depends on
$\delta\boldsymbol{\mu}^T\mathbf{d}$. By adapting the amount of
constraint
tightening $\delta$ to satisfy
\begin{equation}
\label{eq:dualBound} \delta(\boldsymbol{\mu}^k)^T\mathbf{d}\leq \epsilon\ell^*(\bar{x})
\end{equation}
for some $\epsilon>0$ and use this together with the following test
\begin{equation}
\label{eq:perfTestAdapt} D_N^{\delta}(\bar{x},\boldsymbol{\lambda}^k,\boldsymbol{\mu}^k) \geq
P_N(A\bar{x}+Bv_0^k,\mathbf{v}_s^k)+\alpha\ell(\bar{x},v_0^k)
\end{equation}
we get from Lemma~\ref{lem:tighteningRelationsVal} and if
\eqref{eq:dualBound} and 
\eqref{eq:perfTestAdapt} holds that
\ifoneColumn
\begin{align*}
V_N^0(\bar{x})&\geq
D_N^{\delta}(\bar{x},\boldsymbol{\lambda}^k,\boldsymbol{\mu}^k)-\delta(\boldsymbol{\mu}^k)^T\mathbf{d}
 \geq
P_N(A\bar{x}+Bv_0^k,\mathbf{v}_s^k)+\alpha\ell(\bar{x},v_0^k)-\epsilon\ell^*(\bar{x})\\
&\geq V_N^0(A\bar{x}+Bv_0^k)+(\alpha-\epsilon)\ell(\bar{x},v_0^k).
\end{align*}
\else
\begin{align*}
V_N^0(\bar{x})&\geq
D_N^{\delta}(\bar{x},\boldsymbol{\lambda}^k,\boldsymbol{\mu}^k)-\delta(\boldsymbol{\mu}^k)^T\mathbf{d}\\
& \geq
P_N(A\bar{x}+Bv_0^k,\mathbf{v}_s^k)+\alpha\ell(\bar{x},v_0^k)-\epsilon\ell^*(\bar{x})\\
&\geq V_N^0(A\bar{x}+Bv_0^k)+(\alpha-\epsilon)\ell(\bar{x},v_0^k).
\end{align*}
\fi
This is condition \eqref{eq:relDyn}, which guarantees
stability and performance specified by \eqref{eq:perfRes} if
$\alpha>\epsilon$.


\subsection{The stopping condition}

Below we state the stopping condition, whereafter parameter settings
are discussed.
\medskip
\begin{alg}\label{alg_stop_cond} \textbf{Stopping condition}
\hrule
\vspace{2mm}
\noindent {\bf{Input}}: $\bar{x}$\\
Set: $k=0$, $l=0$, $\delta =
\delta_{\rm{init}}$\\
Initialize algorithm \eqref{eq:accGrad1}-\eqref{eq:accGrad4} with:\\
$\boldsymbol{\lambda}^0=\boldsymbol{\lambda}^{-1}=0,
\boldsymbol{\mu}^0=\boldsymbol{\mu}^{-1}=0$ and
$\mathbf{y}^{0}=\mathbf{y}^{-1}=0$.\\
{\bf{Do}}\\
\indent\indent {\bf{If}} $D_{N}^{\delta}(\bar{x},\boldsymbol{\lambda}^k,\boldsymbol{\mu}^k)\geq
P_N(\bar{x},\mathbf{v}^k)-\frac{\epsilon}{l+1}\ell^*(\bar{x})$\\
\indent\indent\indent {\bf{or}} $\delta
\mathbf{d}^T\boldsymbol{\mu}^k > \epsilon\ell^*(\bar{x})$\\
\indent\indent\indent\indent Set
$\delta\leftarrow\delta/2$\phantom{aa}// reduce
constraint tightening\\
\indent\indent\indent\indent Set $l\leftarrow l+1$\\
\indent\indent\indent\indent Set $k=0$ \phantom{aaaa}// reset
step size and iteration counter\\
\indent\indent{\bf{End}}\\
\indent\indent Run $\Delta k$ iterations of \eqref{eq:accGrad1}-\eqref{eq:accGrad4}\\
\indent\indent Set $k\leftarrow k+\Delta k$\\
{\bf{Until}} \begin{tabular}[t]{l}
$D_{N}^{\delta}(\bar{x},\boldsymbol{\lambda}^k,\boldsymbol{\mu}^k)\geq
P_N(A\bar{x}+Bv^k_{0},\mathbf{v}_s^k)+\alpha\ell(\bar{x},v^k_{0})$ {\bf{and}} \\
$\delta\mathbf{d}^T\boldsymbol{\mu}^k \leq \epsilon\ell^*(\bar{x})$
\end{tabular}\\
{\bf{Output}}: $v^k_{0}$
\end{alg}
\medskip
\hrule
\bigskip
In Algorithm~\ref{alg_stop_cond}, four parameters need to be set. 
The first is the
performance parameter $\alpha\in(0,1]$ which guarantees closed loop
performance as specified by \eqref{eq:perfRes}. The larger $\alpha$,
the better performance is guaranteed but a longer control horizon
$N$ will be needed to guarantee the specified performance. The second parameter
is an initial constraint tightening parameter, which we denote by
$\delta_{\rm{init}}\in(0,1]$, from which the constraint tightening
parameter $\delta$ will be adapted (reduced), to satisfy
\eqref{eq:dualBound}. A generic
value that always works is $\delta_{\rm{init}}=0.2$, i.e., $20
\%$ initial constraint tightening. The third parameter
is the relative 
optimality tolerance $\epsilon>0$ where $\epsilon<\alpha$. The
$\epsilon$ must be chosen to satisfy \eqref{eq:alphaBound}. Finally,
$\Delta k$, which is the number of algorithm iterations between every
stopping condition test, should be set to a positive integer,
typically in the range 5 to 20.

Except for the
initial condition $\bar{x}$, Algorithm~\ref{alg_stop_cond} is always 
identically initialized and follows a deterministic scheme. Thus, for
fixed initial condition the same control action is always computed.
This implies that 
Algorithm~\ref{alg_stop_cond} defines a static
feedback control law, which we denote by $\nu_N$. We get the following
closed loop dynamics
\begin{align*}
x_{t+1} &= Ax_{t}+B\nu_N(x_{t}), & x_{0} &= \bar{x}.
\end{align*}
The objective of this section is to present a theorem stating that the
feedback control law function $\nu_N$ satisfies
${\rm{dom}}(\nu_N)\supseteq{\rm{int}}(\mathbb{X}_N^{0})$, where 
\begin{equation}
\label{eq:feasibleRegion}\mathbb{X}_N^{\delta} := \{\bar{x}\in\mathbb{R}^n~|~
V_N^{\delta}(\bar{x})<\infty {\hbox{ and }} Az_{N-1}^*(\bar{x},0)\in{\rm{int}}(\mathcal{X})\}
\end{equation}
which satisfies $\mathbb{X}_N^{\delta_1}\subseteq
\mathbb{X}_N^{\delta_2}$ for $\delta_1 > \delta_2$.
First, however we state the following definition.
\begin{defin}
The constant $\Phi_N$ is the smallest constant such that 
the optimal
solution $\{z_{\tau}^*(\bar{x},0)\}_{\tau=0}^{N-1}$, $\{v_{\tau}^*(\bar{x},0)\}_{\tau=0}^{N-1}$
to \eqref{eq:truncValue} for every $\bar{x}\in \mathbb{X}_N^{0}$
satisfies
\begin{equation*}
\ell^*(z^*_{N-1}(\bar{x},0))\leq\Phi_N\ell(\bar{x},v^*_{0}(\bar{x},0))
\end{equation*}
for the chosen control horizon $N$.
\label{def:ctrbKnown}
\end{defin}

\iflongVersion
The parameter $\Phi_N$ is a measure that compares the first and last
stage costs in the horizon.
In Section~\ref{sec:ctrbVer} a method to compute $\Phi_N$ is presented.
\begin{rem}
In \cite{grimm,grune} an exponential controllability on the stage costs
is assumed, 
i.e., that for $C\geq 1$ and $\sigma\in (0,1)$ the following holds for
$\tau = 0,\ldots,N-1$
\begin{equation*}
\ell^*(z^*_{\tau}(\bar{x},0),v^*_{\tau}(\bar{x},0))\leq C\sigma^{\tau}\ell(\bar{x},v^*_{0}(\bar{x},0)).
\end{equation*}
This implies $\Phi_N \leq C\sigma^{N-1}$.
\label{rem:expCtrb}
\end{rem}
We also need the following lemmas, that are proven in
Appendix\ref{app:finiteIters},
Appendix\ref{app:tighteningRelationsVar} and
Appendix\ref{app:detectInfeas} respectively, to prove the upcoming theorem.
\fi


\begin{lem}
Suppose that $\epsilon>0$ and $\delta\in(0,1]$. For every
$\bar{x}\in\mathbb{X}_N^{\delta}$ we have for some finite $k$ that
\begin{equation}
\label{eq:optCond} D_{N}^{\delta}(\bar{x},\boldsymbol{\lambda}^k,\boldsymbol{\mu}^k)\geq P_N(\bar{x},\mathbf{v}^k)-\epsilon\ell^*(\bar{x}).
\end{equation}
\label{lem:finiteIters}
\end{lem}
\begin{lem}
Suppose that $\epsilon>0$ and $\delta\in(0,1]$. For every
$\bar{x}\in\mathbb{X}_N^{\delta}$ and algorithm iteration $k$ such
that \eqref{eq:optCond} holds we have for $\tau = 
0,\ldots,N-1$ that
\begin{equation*}
\frac{1}{2}\left\|\begin{bmatrix}
\xi^k_{\tau}(\bar{x},\delta)\\
v^k_{\tau}(\bar{x},\delta)
\end{bmatrix}-\begin{bmatrix}
z^*_{\tau}(\bar{x},0)\\
v^*_{\tau}(\bar{x},0)
\end{bmatrix}\right\|_{H}^2 \leq \epsilon\ell^*(\bar{x})+\delta(\boldsymbol{\mu}^k)^T\mathbf{d}
\end{equation*}
where $H = {\rm{blkdiag}}(Q,R)$.
\label{lem:tighteningRelationsVar}
\end{lem}
\begin{lem}
Suppose that $\epsilon>0$ and $\delta\in(0,1]$. For
$\bar{x}\in\mathbb{X}_N^{0}$ but $\bar{x}\notin\mathbb{X}_N^{\delta}$
we have that $\delta (\boldsymbol{\mu}^k)^T\mathbf{d} >
\epsilon\ell^*(\bar{x})$
with finite $k$.
\label{lem:detectInfeas}
\end{lem}

\iflongVersion
We are now ready to state the following theorem, which is proven in
Appendix\ref{app:wellDefStatic}.
\fi

\begin{thm}
Assume that $\epsilon>0$, $\delta_{\rm{init}}\in(0,1]$ and
\begin{equation}
\alpha\leq
1-\epsilon-\kappa(\sqrt{2\epsilon}+\sqrt{\Phi_N})^2(\sqrt{2\epsilon}+1)^2
\label{eq:alphaBound}.
\end{equation}
Then the feedback control law $\nu_N$, defined by
Algorithm~\ref{alg_stop_cond}, satisfies
${\rm{dom}}(\nu_N)\supseteq{\rm{int}}(\mathbb{X}_N^{0})$. Further
\begin{equation}
\label{eq:valueFcnDecreaseThm} V_{N}^{0}(\bar{x})\geq
V_N^0(A\bar{x}+B\nu_N(\bar{x}))+(\alpha-\epsilon)\ell(\bar{x},\nu_N(\bar{x}))
\end{equation}
holds for every $\bar{x}\in\rm{dom}(\nu_N)$.
\label{thm:wellDefStatic}
\end{thm}

\begin{cor}
Suppose that $\alpha\leq 1-\kappa\Phi_N$ and that $\nu_N^*(\bar{x}) =
v_0^*(\bar{x},0)$. Then
\begin{equation*}
V_N^0(\bar{x})\geq V_N^0(A\bar{x}+B\nu_N^*(\bar{x}))+\alpha\ell(\bar{x},\nu_N^*(\bar{x})).
\end{equation*}
holds for every $\bar{x}\in\mathbb{X}_N^{0}$.
\label{cor:optFBperf}
\end{cor}
\iflongVersion
\begin{pf}
For every $\bar{x}\in\mathbb{X}_N^0$ we have
\begin{align*}
V_N^0(\bar{x}) &=
\sum_{\tau=0}^{N-1}\ell(z_{\tau}^*,u_{\tau}^*)
+\ell(Az_{N-1}^*,0)-\ell(Az_{N-1}^*,0)\\
& \geq
V_N^0(A\bar{x}+B\nu_N^*(\bar{x}))+\ell(\bar{x},v_0^*)-\ell(Az_{N-1}^*,0)\\
& \geq
V_N^0(A\bar{x}+B\nu_N^*(\bar{x}))+\ell(\bar{x},v_0^*)-\kappa\ell(z_{N-1}^*,0)\\
& \geq
V_N^0(A\bar{x}+B\nu_N^*(\bar{x}))+(1-\kappa\Phi_N)\ell(\bar{x},v_0^*)
\end{align*}
where the first inequality holds since $Az_{N-1}^*\in\mathcal{X}$ by
construction of $\mathbb{X}_N^0$, the second due to the definition of
$\kappa$ and the third due to the definition of $\Phi_N$.
\end{pf}
\fi

\begin{rem}
By setting $\epsilon=0$ in Theorem~\ref{thm:wellDefStatic} we get
$\alpha\leq 1-\kappa\Phi_N$ as in Corollary~\ref{cor:optFBperf}.
\end{rem}

\subsection{Feasibility, stability and performance}

The following proposition shows one-step feasibility when using the feedback
control law $\nu_N$.
\begin{prp}
Suppose that $\alpha$ satisfies \eqref{eq:alphaBound}.
For every $x_{t}\in{\rm{int}}(\mathbb{X}_N^{0})$ we have that
$x_{t+1}=Ax_{t}+B\nu_N(x_{t})\in\mathcal{X}$.
\end{prp}
\iflongVersion
\begin{pf}
From Theorem~\ref{thm:wellDefStatic} we have that
$x_t\in{\rm{dom}}(\nu_N)$ and from
Algorithm~\ref{alg_stop_cond} we have that 
$P_N(x_{t+1},\mathbf{v}_s^k)<\infty$ which, by definition, implies that
$x_{t+1}\in\mathcal{X}$.
\end{pf}
\fi
The proposition shows that $x_{t+1}$ is feasible if $x_t\in{\rm{dom}}(\nu_N)$. We
define the recursively feasible set as the maximal set such that
\begin{equation*}
\mathbb{X}_{\rm{rf}} = \{x\in\mathcal{X}~|~Ax+B\nu_N(x)\in\mathbb{X}_{\rm{rf}}\}
\end{equation*}
In the following theorem we show that $\mathbb{X}_{\rm{rf}}$ is the region
of attraction and that the control law $\nu_N$ achieves a prespecified
performance as specified by \eqref{eq:perfMetric}.
\begin{thm}
Suppose that $\alpha>\epsilon$ satisfies \eqref{eq:alphaBound}.
Then for every initial condition
$\bar{x}\in\mathbb{X}_{\rm{rf}}$ we have that $\|x_{t}\|\to 0$ as
$t\to\infty$ and that the
closed loop performance satisfies
\begin{equation}
(\alpha-\epsilon) V_{\infty,\nu_N}(\bar{x})\leq V_{\infty}(\bar{x}).
\label{eq:clPerf}
\end{equation}
Further, $\mathbb{X}_{\rm{rf}}$ is the region of attraction.
\label{thm:fixedStab}
\end{thm}
\iflongVersion
\begin{pf}
From the definition of $\mathbb{X}_{\rm{rf}}$ we know that
$\bar{x}=x_{0}\in\mathbb{X}_{\rm{rf}}$ implies
$x_{t}\in\mathbb{X}_{\rm{rf}}$ for all $t\in\mathbb{N}_{\geq 0}$.
Since, by construction, $\mathbb{X}_{\rm{rf}}\subseteq
{\rm{int}}(\mathbb{X}_N^0)\subseteq{\rm{dom}}(\nu_N)$ we have from
Theorem~\ref{thm:wellDefStatic} that
\eqref{eq:valueFcnDecreaseThm} holds for all $x_{t}$,
$t\in\mathbb{N}_{\geq 0}$.
In \cite[Proposition 2.2]{gruneRantzer} it was shown, using telescope
summation, that \eqref{eq:valueFcnDecreaseThm} implies \eqref{eq:clPerf}.
Further, 
since the stage cost $\ell$ satisfies \cite[Assumption 5.1]{grune}
we get from \cite[Theorem 5.2]{grune} that $\|x_{t}\|\to 0$ as
$t\to\infty$.

What is left to show is that $\mathbb{X}_{\rm{rf}}$ is the region of
attraction. Denote by $\mathbb{X}_{\rm{roa}}$ the region of
attraction using $\nu_N$. We have above shown that
$\mathbb{X}_{\rm{rf}}\subseteq\mathbb{X}_{\rm{roa}}$. We next show that
$\mathbb{X}_{\rm{roa}}\subseteq\mathbb{X}_{\rm{rf}}$ by a 
contradiction argument to conclude that
$\mathbb{X}_{\rm{rf}}=\mathbb{X}_{\rm{roa}}$. Assume that
there exist $\bar{x}\in\mathbb{X}_{\rm{roa}}$ such that
$\bar{x}\notin\mathbb{X}_{\rm{rf}}$. If $\bar{x}\in\mathbb{X}_{\rm{roa}}$ the
closed loop state sequence $\{x_t\}_{t=0}^\infty$ is feasible in every
step (and
converges to the origin) and consequently 
$\{Ax_t+B\nu_N(x_t)\}_{t=0}^\infty$ is feasible in every step. This is exactly the
requirement to have $\bar{x}\in\mathbb{X}_{\rm{rf}}$, which is a
contradiction. Thus
$\mathbb{X}_{\rm{rf}}\subseteq\mathbb{X}_{\rm{roa}}\subseteq\mathbb{X}_{\rm{rf}}$
which implies that
$\mathbb{X}_{\rm{rf}}=\mathbb{X}_{\rm{roa}}$.

This completes the proof.
\end{pf}
\fi

\iflongVersion
To guarantee a priori that the control law $\nu_N$ achieves the
performance \eqref{eq:clPerf} specified by $\alpha$, we need to find a
control horizon $N$ such that the corresponding controllability
parameter $\Phi_N$ satisfies
\eqref{eq:alphaBound}. This requires
the computation of controllability parameter $\Phi_N$ which is the
topic of the next section.
\else
\fi

\iflongVersion
\section{Offline controllability verification}
\label{sec:ctrbVer}
The stability and performance results in Theorem~\ref{thm:fixedStab}
rely on Definition~\ref{def:ctrbKnown}. For the results to 
be practically meaningful it must be possible to compute $\Phi_N$ in
Definition~\ref{def:ctrbKnown}. In this section we will show that
this can be done by solving a mixed integer linear program (MILP).
For desired performance specified by $\alpha$, we get a requirement
on the controllability parameter through \eqref{eq:alphaBound} for
Theorem~\ref{thm:wellDefStatic} and Theorem~\ref{thm:fixedStab} to hold. We
denote by $\Phi_\alpha$ the largest controllability parameter such
that Theorem~\ref{thm:wellDefStatic} and Theorem~\ref{thm:fixedStab}
holds for the specified $\alpha$. This parameter is the one that gives
equality in \eqref{eq:alphaBound}, i.e., satisfies
\begin{equation}
\alpha= 1-\epsilon-\kappa
(\sqrt{2\epsilon}+\sqrt{\Phi_\alpha})^2(\sqrt{2\epsilon}+1)^2
\label{eq:alphaBoundExact}
\end{equation}
for the desired performance $\alpha$ and optimality tolerance
$\epsilon$.
The parameters $\alpha$ and $\epsilon$ must be chosen such
that $\Phi_{\alpha}>0$. The objective is to find a control horizon $N$
such that the corresponding 
controllability parameter $\Phi_N$ satisfies $\Phi_N\leq\Phi_\alpha$.
First we show that for long enough control horizon $N$ there
exist a $\Phi_N\leq\Phi_\alpha$.
\begin{lem}
Assume that $\alpha$ and $\epsilon$ are chosen such that
$\Phi_\alpha>0$ where $\Phi_\alpha$ is implicitly defined in
\eqref{eq:alphaBoundExact}. Then there 
exists control horizon $N$ and corresponding controllability
parameter $\Phi_N\leq\Phi_\alpha$.
\end{lem}
\begin{pf}
Since $\mathbb{X}_{\rm{rf}}$ is the region of
attraction we have $\mathbb{X}_{\rm{rf}}\subseteq\mathbb{X}_{\infty}$. This
in turn implies that
\eqref{eq:truncValueCompact} is feasible for every control horizon
$N\in\mathbb{N}_{\geq 1}$
due to the absence of terminal constraints.
We have
\ifoneColumn
\begin{align*}
V_N(\bar{x}) &= \sum_{\tau = 0}^{N-2}\ell(z^*_{\tau},v^*_{\tau})+\ell(z^*_{N-1},v^*_{N-1})
\geq V_{N-1}(\bar{x})+\ell(z^*_{N-1},v^*_{N-1}).
\end{align*}
\else
\begin{align*}
V_N(\bar{x}) &= \sum_{\tau = 0}^{N-2}\ell(z^*_{\tau},v^*_{\tau})+\ell(z^*_{N-1},v^*_{N-1})\\
&\geq V_{N-1}(\bar{x})+\ell(z^*_{N-1},v^*_{N-1}).
\end{align*}
\fi
Since the pair $(A,B)$ is assumed controllable and since
\eqref{eq:truncValueCompact} has neither terminal constraints nor
terminal cost we have for some finite $M$ that
$M\geq V_\infty(\bar{x})\geq V_N(\bar{x})\geq V_{N-1}(\bar{x})$. Thus the sequence
$\{V_N(\bar{x})\}_{N=0}^{\infty}$ is a bounded monotonic increasing sequence
which is well known to be convergent. Thus, for $N\geq \bar{N}$ where
$\bar{N}$ is large enough the
difference $V_N(\bar{x})-V_{N-1}(\bar{x})$ is arbitrarily small. Especially
$\ell(z^*_{N-1},v^*_{N-1})=\ell^*(z^*_{N-1})\leq V_N(\bar{x})-V_{N-1}(\bar{x})\leq
\Phi_\alpha\ell(\bar{x},v^*_{0})$ since $\Phi_\alpha>0$. That is,
for long enough control horizon
$N\geq \bar{N}$, $\Phi_N\leq\Phi_\alpha$. This
completes the proof.
\end{pf}

The preceding Lemma shows that there exists a control horizon $N$
such that $\Phi_N\leq\Phi_\alpha$ if $\Phi_\alpha>0$ for the chosen
performance $\alpha$ and tolerance $\epsilon$.
The choice of performance parameter $\alpha$ gives requirements on how
$\epsilon$ can be chosen to give $\Phi_\alpha>0$. Larger $\epsilon$ requires smaller
$\Phi_{\alpha}$ to satisfy \eqref{eq:alphaBoundExact} which in turn
requires longer control horizons $N$ since $\Phi_N$ must satisfy
$\Phi_N\leq\Phi_{\alpha}$. In the
following section we address the problem of how to compute the
control horizon $N$ and corresponding $\Phi_N$ such that the desired
performance specified by $\alpha$ can be guaranteed.

\subsection{Exact verification of controllability parameter}

In the following proposition we introduce an optimization problem that
tests if the controllability parameter $\Phi_N$ corresponding to
control horizon $N$ satisfies $\Phi_N\leq\Phi_\alpha$ for the desired
performance specified by $\alpha$.
Before we state the
proposition, the following matrices are introduced
\begin{align*}
T &= {\rm{blkdiag}}(0,\ldots,0,-Q,\Phi_\alpha R,0,\ldots,0,-R)\\
S &= {\rm{blkdiag}}(0,\ldots,0,I,0,\ldots,0)
\end{align*}
where $Q$ and $R$ are
the cost matrices for states and inputs and $\Phi_\alpha$ is the
required controllability parameter for the chosen $\alpha$.
Recalling the partitioning \eqref{eq:chiDef} of $\mathbf{y}$
implies that 
\begin{align*}
\mathbf{y}^TT\mathbf{y} &=
\Phi_{\alpha}v_{0}^TRv_0-z_{N-1}^TQz_{N-1}-v_{N-1}^TRv_{N-1}\\
S\mathbf{y} &= z_{N-1}
\end{align*}
\begin{prp}
Assume that $\Phi_{\alpha}>0$ satisfies \eqref{eq:alphaBoundExact} for
the chosen performance parameter $\alpha$ and optimality tolerance
$\epsilon$. Further assume that the control horizon $N$ is such that 
\begin{align}
\label{eq:bilevel} 0=\min_{\bar{x}} ~&\frac{1}{2}\left(\Phi_\alpha\bar{x}^TQ\bar{x}+\mathbf{y}^TT\mathbf{y}\right)\\
\nonumber {\rm{s.t.}}~& \bar{x}\in\mathbb{X}_N^0\\
\nonumber ~& \mathbf{y} =\arg\min V_{N}^{0}(\bar{x})
\end{align}
then $\Phi_N\leq \Phi_\alpha$.
\label{prp:bilevel}
\end{prp}
\begin{pf}
First we note that $\bar{x}=0$ gives $\mathbf{y}=0$ and
$\Phi_\alpha\bar{x}^TQ\bar{x}+\mathbf{y}^TT\mathbf{y}=0$,
i.e., we have that 0 is always a feasible solution.
Further, \eqref{eq:bilevel} implies for every $\bar{x}\in\mathbb{X}_N^0$ that
\ifoneColumn
\begin{align*}
0&\leq\Phi_\alpha\bar{x}^TQ\bar{x}+\mathbf{y}^TT\mathbf{y}
=\Phi_\alpha\ell(\bar{x},v^*_{0})-\ell(z^*_{N-1},v^*_{N-1})
=\Phi_\alpha\ell(\bar{x},v^*_{0})-\ell^*(z^*_{N-1})
\end{align*}
\else
\begin{align*}
0&\leq\Phi_\alpha\bar{x}^TQ\bar{x}+\mathbf{y}^TT\mathbf{y}
=\Phi_\alpha\ell(\bar{x},v^*_{0})-\ell(z^*_{N-1},v^*_{N-1})\\
&=\Phi_\alpha\ell(\bar{x},v^*_{0})-\ell^*(z^*_{N-1})
\end{align*}
\fi
since $v^*_{N-1}=0$. This is exactly the condition in
Definition~\ref{def:ctrbKnown}. Since $\Phi_N$ is the smallest such
constant, we have 
$\Phi_N\leq\Phi_\alpha$ for the chosen 
control horizon $N$ and desired performance $\alpha$ and 
optimality tolerance $\epsilon$.
\end{pf}

The optimization problem \eqref{eq:bilevel} is a bilevel
optimization problem with indefinite quadratic cost (see \cite{colson}
for a survey on bilevel optimization). Such problems are
in general NP-hard to solve.
The problem can, however, be rewritten as an equivalent MILP
as shown in the following proposition which is a straightforward
application of \cite[Theorem 2]{JonesBilevel}.
\begin{prp}
Assume that $\Phi_{\alpha}$ satisfies \eqref{eq:alphaBoundExact} for
the chosen performance parameter $\alpha$ and optimality tolerance
$\epsilon$. If the control horizon $N$ is such that the following holds
\ifoneColumn
\begin{align}
\label{eq:MILP} 0=\min ~&-\frac{1}{2}\left(d_x^T\mu^{U1}+d_x^T\mu^{U2}+\mathbf{d}^T\mu^{UL1}\right)\\
\nonumber{\rm{s.t.}}~&\beta_i^L\in\{0,1\}~,~\beta_i^{U1}\in\{0,1\}~,~\beta_i^{U2}\in\{0,1\}\\
\nonumber~& {\hbox{\small{Upper level}}}\\
\nonumber~&\left\lfloor\begin{array}{l}
{\hbox{\small{Primal and dual feasibility}}}\\
\left\lfloor\begin{array}{l}
C_x\bar{x}-d_x-s^x=0\\
 s^x\leq 0~,~\mu^{U1}\geq 0\\
C_xAS\mathbf{y}-d_x-s^{z} = 0\\
 s^{z}\leq 0~,~\mu^{U2}\geq 0\\
\end{array}\right.\\
{\hbox{\small{Stationarity}}}\\
\left\lfloor\begin{array}{l}
\Phi_\alpha Q\bar{x}+(C_x)^T\mu^{U1}-\mathbf{b}^T\lambda^{UL2}=0\\
T\mathbf{y}+H^T\lambda^{UL1}+\mathbf{A}^T\lambda^{UL2}+\mathbf{C}^T\mu^{UL1}
+(C_xAS)^T\mu^{U2}=0\\
 \mathbf{A}\lambda^{UL1} = 0\\
 \mathbf{C}\lambda^{UL1}-\mu^{UL2}=0\\
\end{array}\right.\\
{\hbox{\small{Complementarity}}}\\
\left\lfloor\begin{array}{l}
\beta_i^L = 1\Rightarrow \mu_i^{UL2}=0~,~ \beta_i^L = 0\Rightarrow \mu_i^{UL1}=0\\
\beta_i^{U1} = 1\Rightarrow s^x_i=0~,~ \beta_i^{U1} = 0\Rightarrow
\mu_i^{U1}=0\\
\beta_i^{U2} = 1\Rightarrow s^z_i=0~,~ \beta_i^{U2} = 0\Rightarrow \mu_i^{U2}=0\\
\end{array}\right.\\
\end{array}\right.\\
\nonumber~& {\hbox{\small{Lower level}}}\\
\nonumber~& \left\lfloor\begin{array}{l}
{\hbox{\small{Primal and dual feasibility}}}\\
\left\lfloor\begin{array}{l}
 \mathbf{A}\mathbf{y}-\mathbf{b}\bar{x}=0\\
\mathbf{C}\mathbf{y}-\mathbf{d}-s =0\\
 s\leq 0~,~\mu^L \geq 0\\
\end{array}\right.\\
{\hbox{\small{Stationarity}}}\\
\left\lfloor\begin{array}{l}
H\mathbf{y}+\mathbf{A}^T\lambda^L+\mathbf{C}^T\mu^{L} =0\\
\end{array}\right.\\
{\hbox{\small{Complementarity}}}\\
\left\lfloor\begin{array}{l}
\beta_i^L = 1\Rightarrow s_i=0~,~ \beta_i^L = 0\Rightarrow \mu_i^L=0\\
\end{array}\right.\\
\end{array}\right.
\end{align}
\else
\begin{align}
\label{eq:MILP} 0=\min ~&-\frac{1}{2}\left(d_x^T\mu^{U1}+d_x^T\mu^{U2}+\mathbf{d}^T\mu^{UL1}\right)\\
\nonumber{\rm{s.t.}}~&\beta_i^L\in\{0,1\}~,~\beta_i^{U1}\in\{0,1\}~,~\beta_i^{U2}\in\{0,1\}\\
\nonumber~& {\hbox{\small{Upper level}}}\\
\nonumber~&\left\lfloor\begin{array}{l}
{\hbox{\small{Primal and dual feasibility}}}\\
\left\lfloor\begin{array}{l}
C_x\bar{x}-d_x-s^x=0\\
 s^x\leq 0~,~\mu^{U1}\geq 0\\
C_xAS\mathbf{y}-d_x-s^{z} = 0\\
 s^{z}\leq 0~,~\mu^{U2}\geq 0\\
\end{array}\right.\\
{\hbox{\small{Stationarity}}}\\
\left\lfloor\begin{array}{l}
\Phi_\alpha Q\bar{x}+(C_x)^T\mu^{U1}-\mathbf{b}^T\lambda^{UL2}=0\\
T\mathbf{y}+H^T\lambda^{UL1}+\mathbf{A}^T\lambda^{UL2}+\mathbf{C}^T\mu^{UL1}\\
\qquad\qquad\qquad\qquad\qquad+(C_xAS)^T\mu^{U2}=0\\
 \mathbf{A}\lambda^{UL1} = 0\\
 \mathbf{C}\lambda^{UL1}-\mu^{UL2}=0\\
\end{array}\right.\\
{\hbox{\small{Complementarity}}}\\
\left\lfloor\begin{array}{l}
\beta_i^L = 1\Rightarrow \mu_i^{UL2}=0~,~ \beta_i^L = 0\Rightarrow \mu_i^{UL1}=0\\
\beta_i^{U1} = 1\Rightarrow s^x_i=0~,~ \beta_i^{U1} = 0\Rightarrow
\mu_i^{U1}=0\\
\beta_i^{U2} = 1\Rightarrow s^z_i=0~,~ \beta_i^{U2} = 0\Rightarrow \mu_i^{U2}=0\\
\end{array}\right.\\
\end{array}\right.\\
\nonumber~& {\hbox{\small{Lower level}}}\\
\nonumber~& \left\lfloor\begin{array}{l}
{\hbox{\small{Primal and dual feasibility}}}\\
\left\lfloor\begin{array}{l}
 \mathbf{A}\mathbf{y}-\mathbf{b}\bar{x}=0\\
\mathbf{C}\mathbf{y}-\mathbf{d}-s =0\\
 s\leq 0~,~\mu^L \geq 0\\
\end{array}\right.\\
{\hbox{\small{Stationarity}}}\\
\left\lfloor\begin{array}{l}
H\mathbf{y}+\mathbf{A}^T\lambda^L+\mathbf{C}^T\mu^{L} =0\\
\end{array}\right.\\
{\hbox{\small{Complementarity}}}\\
\left\lfloor\begin{array}{l}
\beta_i^L = 1\Rightarrow s_i=0~,~ \beta_i^L = 0\Rightarrow \mu_i^L=0\\
\end{array}\right.\\
\end{array}\right.
\end{align}
\fi
where all $\beta, \mu, \lambda, s$ and $\bar{x}, \mathbf{y}$
are decision variables,
then $\Phi_\alpha\geq\Phi_N$. 
\label{prp:MILP}
\end{prp}
\begin{pf}
The set $\mathbb{X}_N^0$ can equivalently be written as
\ifoneColumn
\begin{align}
\label{eq:altDefXN}\mathbb{X}_N^0 = \{x\in\mathbb{R}^n~|~&
\mathbf{A}\mathbf{y}^*(x,0)=\mathbf{b}x,
\mathbf{C}\mathbf{y}^*(x,0)\leq\mathbf{d},
C_xAS\mathbf{y}^*(x,0)\leq d_x, C_xx\leq d_x\}.
\end{align}
\else
\begin{align}
\nonumber\mathbb{X}_N^0 = \{x\in\mathbb{R}^n~|~&
\mathbf{A}\mathbf{y}^*(x,0)=\mathbf{b}x,
\mathbf{C}\mathbf{y}^*(x,0)\leq\mathbf{d},\\
\label{eq:altDefXN} &C_xAS\mathbf{y}^*(x,0)\leq d_x, C_xx\leq d_x\}.
\end{align}
\fi
We express the set
$\mathbb{X}_N^0$ in \eqref{eq:bilevel} using
\eqref{eq:altDefXN}. The equivalence
between the optimization problems \eqref{eq:MILP} and \eqref{eq:bilevel} is
established in \cite[Theorem 2]{JonesBilevel}. The remaining parts of
the proposition follow by applying Proposition~\ref{prp:bilevel}.
\end{pf}

The transformation from \eqref{eq:bilevel} to \eqref{eq:MILP} is done
by expressing the lower level optimization problem in
\eqref{eq:bilevel} by its sufficient and
necessary KKT conditions to get a single level indefinite quadratic
program with
complementarity constraints. The resulting indefinite quadratic
program with complementarity constraints can in turn be cast as a
MILP to get \eqref{eq:MILP}.
\begin{rem}
Although MILP problems are NP-hard, there are efficient solvers
available such as CPLEX and GUROBI.
There are also solvers available for solving the bilevel optimization
problem \eqref{eq:bilevel} directly, e.g., the function
{\it{solvebilevel}} in YALMIP, \cite{YALMIP}. 
\end{rem}

If the chosen control horizon $N$ is not long enough for $\Phi_N\leq
\Phi_\alpha$, different heuristics
can be used to choose a new longer horizon to be verified. One
heuristic is to assume exponential controllability as in
Remark~\ref{rem:expCtrb}, i.e., that there exist constants $C\geq 1$ and
$\sigma\in(0,1)$ such that
\begin{equation}
\label{eq:expCtrb} C\sigma^\tau\ell(\bar{x},v_0^k)\geq\ell(z_\tau^k,v_\tau^k)
\end{equation}
for all $\tau = 0,\ldots,N-1$. The $C$ and $\sigma$-parameters
should be determined using the optimal solution
$\mathbf{y}$ to \eqref{eq:truncValueCompact} for the $x$ that
minimized \eqref{eq:MILP} in the previous test. Under the assumption
that \eqref{eq:expCtrb} holds as $N$ increases, a new
guess on the control horizon $N$ can be computed by finding the
smallest $N$ such that $C\sigma^{N-1}\leq \Phi_\alpha$.

\subsection{Controllability parameter estimation}

The test in Proposition~\ref{prp:MILP} verifies if the control horizon $N$ is
long enough for the controllability assumption to hold for
the required controllability parameter $\Phi_\alpha$. Thus, an
initial guess on the control horizon is needed. A guaranteed lower bound can
easily be computed by solving \eqref{eq:truncValueCompact} for a
variety of initial conditions $\bar{x}$ and compute the worst
controllability parameter, denoted by $\hat{\Phi}_N$, for these sample
points. If the estimated controllability
parameter $\hat{\Phi}_N\geq\Phi_\alpha$, we know that the control
horizon need to be increased for \eqref{eq:MILP} to hold. If instead 
$\hat{\Phi}_N\leq\Phi_\alpha$ the control horizon $N$ might serve as a
good initial guess to be verified by \eqref{eq:MILP}.

\begin{rem}
For large systems, \eqref{eq:MILP} may be too complex to verify
the desired performance. In such cases, the heuristic method mentioned
above can be 
used in conjunction with an adaptive horizon scheme. The adaptive
scheme keeps the horizon fixed for all time-steps until the
controllability assumption does not hold. Then, the control
horizon is increased to satisfy the assumption and kept at the new
level until the controllability assumption does not hold again.
Eventually the control horizon will be large enough for
$\Phi_N\leq\Phi_\alpha$ and the horizon need not be increased again.
\end{rem}
\else
\fi





\section{Numerical example}
\label{sec:numerical}
We evaluate the efficiency of the proposed distributed feedback
control law $\nu_N$ 
by applying it to a randomly generated dynamical
system with sparsity structure that is specified in \cite[Supplement
A.1]{gisPhD}. The random 
dynamics matrix is scaled such that the magnitude of the largest
eigenvalue is 1.1.
The system has 3 subsystems with 5
states and 1 input each. All
state variables are upper and
lower bounded by random numbers in the intervals $[0.5, 1.5]$ and
$[-0.15, -0.05]$ respectively and all input variables are upper and
lower bounded by random numbers in the intervals $[0.5, 1.5]$ and
$[-0.5, -1.5]$ respectively. The stage cost is chosen to be
\begin{equation*}
\ell_i(x_i,u_i) = x_i^Tx_i+u_i^Tu_i
\end{equation*}
for $i=1,2,3$.
The suboptimality parameter is chosen $\alpha =
0.01$. According to Theorem~\ref{thm:wellDefStatic}, to quantify the control horizon 
$N(\alpha)$, the optimality
tolerance $\epsilon$ must be chosen and $\kappa$ 
computed, where 
$\kappa$ is the smallest constant such that
$\kappa Q\succeq A^TQA$. We get $\kappa = 1.22$ and
choose $\epsilon = 0.005$. Using \eqref{eq:alphaBound}, we get
$\Phi_{N(0.01)} \leq 0.51$. Verification by solving the MILP in
\iflongVersion
\eqref{eq:MILP}
\fi
gives that the smallest control horizon $N(0.01)$ that satisfies
$\Phi_{N(0.01)} \leq 0.51$ is $N(0.01) = 6$.

Table~\ref{tab:results} presents the results. The first
column specifies the stopping condition used, ``stop. cond.'' for the
stopping condition presented in Algorithm~\ref{alg_stop_cond} and ``opt.
cond.'' for a optimality conditions.
The second column specifies the duality gap tolerance $\epsilon$ and the third
column specifies the initial constraint tightening
$\delta_{\rm{init}}$ for the stopping condition and the relative
accuracy requirement for the constraints when using optimality conditions.

\begin{table}
\centering
\caption{Experimental results for different performance requirements
  $\alpha$ and different initial constraint tightenings
  $\delta_{\rm{init}}$.}
\begin{tabular}{|c|c|c|c|c|c|}
\multicolumn{6}{c}{Algorithm comparison, $\alpha=0.01$, $N=6$}\\
\hline
condition & $\epsilon$ & $\delta_{\rm{init}}$ & avg. $\#$ iters & max $\#$ iters & avg.
$\delta$  \\
\hline
stop. cond. &  0.005 & 0.001 & 288.3 & 506 & 0.001 \\
stop. cond. &  0.005 & 0.01 & 151.5 & 260 & 0.01 \\
stop. cond. &  0.005 & 0.05 & 73.7 & 237 & 0.049  \\
stop. cond. &  0.005 & 0.1 & 70.7 & 236 & 0.057  \\
stop. cond. &  0.005 & 0.2 & 72.8 & 236 & 0.060 \\
stop. cond. &  0.005 & 0.5 & 69.2 & 234 & 0.076 \\
\hline
opt. cond.  &  0.005 & 0.001 & 324.5 & 506 & 0.001 \\
opt. cond.  &  0.005 & 0.01 & 171.5 & 260 & 0.01 \\
\hline
\end{tabular}\\
\vspace{3mm}
\label{tab:results}
\end{table}

Columns four, five and six contain the simulation results. The results are
obtained by simulating the system with 1000 randomly chosen initial
conditions that are drawn from a uniform distribution on $\mathcal{X}$.
Column four and five contain the mean and max numbers of iterations needed and
column six presents the average constraint tightening $\delta$
used at termination of Algorithm~\ref{alg_stop_cond}.

We see that the adaptive constraint tightening approach gives 
considerably less iterations for a larger initial tightening. However,
for more than $10 \%$ initial constraint tightening ($\delta_{\rm{init}} = 0.1$),
the number of iterations is not significantly affected. It is
remarkable to note that $50 \%$ initial constraint tightening ($\delta_{\rm{init}} =
0.5$) is as efficient
as, e.g., $5 \%$ ($\delta_{\rm{init}} = 0.05$) considering that more 
reductions in the constraint tightening need to be performed. This
indicates early detection of
infeasibility. We also note that for a suitable choice of initial
constraint tightening, the average number of iterations is
reduced significantly.


\section{Conclusions}\label{sec:conclusions}

We have equipped the duality-based distributed optimization algorithm
in \cite{gisAutomatica},
when used in a DMPC context, with a stopping condition that
guarantees feasibility of the optimization problem and stability and a
prespecified performance of the closed-loop system. 
A numerical example is provided that shows that the stopping
condition can reduce significantly the number of
iterations needed to achieve these properties.

\section{Acknowledgments}\label{sec:acknowledgements}
The authors were supported by the Swedish Research
Council through the Linnaeus center LCCC and the eLLIIT Excellence
Center at Lund University.

\bibliographystyle{IEEEtran}

\bibliography{report}

\iflongVersion
\appendices

\subsection{Proof for Lemma~\ref{lem:finiteIters}}\label{app:finiteIters}

We divide the proof into two parts, the first for $\bar{x}=0$ and the second
for $\bar{x}\neq 0$. For $\bar{x}=0$ we have at iteration $k=0$ that
$\mathbf{y}^0=0$ which is the optimal solution. Hence
\eqref{eq:optCond} holds for $k=0$ since all terms are 0 and
$0=A\xi^0_{N-1}\in\mathcal{X}$.

Next, we show the result for $\bar{x}\neq 0$. Whenever
\eqref{eq:truncValueCompact} is feasible we have
convergence in primal variables \cite[Theorem 1]{gisAutomatica}. This together with
the linear relation through which $\xi$ is defined \eqref{eq:xiDef} gives
$\xi^k_{\tau}\to z^*_{\tau}$ for $\tau=0,\ldots,N-1$ as
$k\to\infty$. We have $z^*_{\tau}\in(1-\delta)\mathcal{X}$
and since $(1-\delta)\mathcal{X}\subset\mathcal{X}$ for every
$\delta\in(0,1]$
this implies that there exists finite $k_0^x$ such that
$\xi^k_{\tau}\in\mathcal{X}$ for all $k\geq k_0^x$. 
Equivalent convergence reasoning holds for $v^k_{\tau}$.
Together this 
implies that there exists finite $k_0^P$ such that
$P_N(\bar{x},\mathbf{v}^k)<\infty$ and that $P_N(\bar{x},\mathbf{v}^k)\to
V_{N}^{\delta}(\bar{x})$ for all $k\geq k_0^P$. Together with convergence
in dual function value \cite[Theorem 1]{gisAutomatica} gives that
\begin{equation*}
D_{N}^{\delta}(\bar{x},\boldsymbol{\lambda}^k,\boldsymbol{\mu}^k)\geq
P_N(\bar{x},\mathbf{v}^k)-\epsilon\ell^*(\bar{x})
\end{equation*}
holds with finite $k$ since $\ell^*(\bar{x})>0$ and $\epsilon>0$. This
concludes the proof.\hfill$\Box$

\subsection{Proof for
  Lemma~\ref{lem:tighteningRelationsVar}}\label{app:tighteningRelationsVar}
We introduce $\mathbf{y}^k = [(\boldsymbol{\xi}^k(\bar{x},\delta))^T
(\mathbf{v}^k(\bar{x},\delta))^T]^T$, where $\boldsymbol{\xi}^k(\bar{x},\delta)$
and $\mathbf{v}^k(\bar{x},\delta)$ satisfies the dynamic equations
\eqref{eq:xiDef}. Whenever \eqref{eq:optCond} holds we have that
$\xi_{\tau}^k(\bar{x},\delta)\in\mathcal{X}$ and
$v_{\tau}^k(\bar{x},\delta)\in\mathcal{U}$ for $\tau = 0,\ldots,N-1$.
We also introduce
 $\mathbf{y}^* = [(\boldsymbol{z}^*(\bar{x},0))^T
(\mathbf{v}^*(\bar{x},0))^T]^T$. This implies
\ifoneColumn
\begin{align*}
\frac{1}{2}(\mathbf{y}^k-\mathbf{y}^*)^T\mathbf{H}(\mathbf{y}^k-\mathbf{y}^*)&=
\frac{1}{2}(\mathbf{y}^k)^T\mathbf{H}\mathbf{y}^k
-\frac{1}{2}(\mathbf{y}^*)^T\mathbf{H}\mathbf{y}^* -\langle
H\mathbf{y}^*,\mathbf{y}^k-\mathbf{y}^*\rangle\\
&\leq P_N(\bar{x},\mathbf{v}^k)-V_N^0(\bar{x})
\leq
D_N^{\delta}(\bar{x},\boldsymbol{\lambda}^k,\boldsymbol{\mu}^k)+\epsilon\ell^*(\bar{x})-V_N^0(\bar{x})\\
&\leq \delta(\boldsymbol{\mu}^k)^T\mathbf{d}+\epsilon\ell^*(\bar{x})
\end{align*}
\else
\begin{align*}
&\frac{1}{2}(\mathbf{y}^k-\mathbf{y}^*)^T\mathbf{H}(\mathbf{y}^k-\mathbf{y}^*)=\\
&=\frac{1}{2}(\mathbf{y}^k)^T\mathbf{H}\mathbf{y}^k
-\frac{1}{2}(\mathbf{y}^*)^T\mathbf{H}\mathbf{y}^* -\langle
H\mathbf{y}^*,\mathbf{y}^k-\mathbf{y}^*\rangle\\
&\leq P_N(\bar{x},\mathbf{v}^k)-V_N^0(\bar{x})
\leq
D_N^{\delta}(\bar{x},\boldsymbol{\lambda}^k,\boldsymbol{\mu}^k)+\epsilon\ell^*(\bar{x})-V_N^0(\bar{x})\\
&\leq \delta(\boldsymbol{\mu}^k)^T\mathbf{d}+\epsilon\ell^*(\bar{x})
\end{align*}
\fi
where the first inequality comes from the first order optimality condition
\cite[Theorem 2.2.5]{NesterovLectures} and by definition of 
$V_N^0$ and $P_N$. The second inequality is due to \eqref{eq:optCond}
and the last inequality follows from Lemma~\ref{lem:tighteningRelationsVal}.
Further, since $\mathbf{H}={\rm{blkdiag}}(Q,\ldots,Q,R,\ldots,R)$ we have 
for $\tau=0,\ldots,N-1$ that
\ifoneColumn
\begin{align*}
\frac{1}{2}\left\|\begin{bmatrix}
\xi^k_{\tau}(\bar{x},\delta)\\
v^k_{\tau}(\bar{x},\delta)
\end{bmatrix}-\begin{bmatrix}
z^*_{\tau}(\bar{x},0)\\
v^*_{\tau}(\bar{x},0)
\end{bmatrix}\right\|_{H}^2&
\leq\frac{1}{2}(\mathbf{y}^k-\mathbf{y}^*)^T\mathbf{H}
(\mathbf{y}^k-\mathbf{y}^*)
\leq \delta(\boldsymbol{\mu}^k)^T\mathbf{d}+\epsilon\ell^*(\bar{x})
\end{align*}
\else
\begin{align*}
\frac{1}{2}\left\|\begin{bmatrix}
\xi^k_{\tau}(\bar{x},\delta)\\
v^k_{\tau}(\bar{x},\delta)
\end{bmatrix}-\begin{bmatrix}
z^*_{\tau}(\bar{x},0)\\
v^*_{\tau}(\bar{x},0)
\end{bmatrix}\right\|_{H}^2&
\leq\frac{1}{2}(\mathbf{y}^k-\mathbf{y}^*)^T\mathbf{H}
(\mathbf{y}^k-\mathbf{y}^*)\\
&\leq \delta(\boldsymbol{\mu}^k)^T\mathbf{d}+\epsilon\ell^*(\bar{x})
\end{align*}
\fi
where $H = {\rm{blkdiag}}(Q,R)$, whenever \eqref{eq:optCond} holds.
This completes the proof.\hfill$\Box$


\subsection{Proof for
  Lemma~\ref{lem:detectInfeas}}\label{app:detectInfeas}
Since $x\in\mathbb{X}_N^0$ but $x\notin\mathbb{X}_N^\delta$ we have
that $V_{N}^{0}(\bar{x})<\infty$ and $V_{N}^{\delta}(\bar{x})=\infty$.
Further, from the
strong theorem of alternatives \cite[Section 5.8.2]{Boyd2004} we know that
since $V_N^{\delta}(\bar{x})=\infty$ for the current constraint
tightening $\delta$ the
dual problem is unbounded. Hence there exist $\boldsymbol{\lambda}_f$,
$\boldsymbol{\mu}_f$ such that
\begin{align}
\label{eq:dualUnbounded}\delta\boldsymbol{\mu}_f^T\mathbf{d}&\geq
D_{N}^{\delta}(\bar{x},\boldsymbol{\lambda}_f,\boldsymbol{\mu}_f)-V_{N}^{0}(\bar{x})\geq 2\epsilon\ell^*(\bar{x})
\end{align}
where Lemma~\ref{lem:tighteningRelationsVal} is used in the first
inequality. Further, the convergence rate in \cite[Theorem
4.4]{BecTab_FISTA:2009} for algorithm
\eqref{eq:accGrad1}-\eqref{eq:accGrad4} is
\begin{equation*}
D_{N}^{\delta}(\bar{x},\boldsymbol{\lambda}^*,\boldsymbol{\mu}^*)-D_{N}^{\delta}(\bar{x},\boldsymbol{\lambda}^k,\boldsymbol{\mu}^k)\leq
\frac{2L}{(k+1)^2}\left\|\begin{bmatrix}
\boldsymbol{\lambda}^*\\
\boldsymbol{\mu}^*
\end{bmatrix}-\begin{bmatrix}
\boldsymbol{\lambda}^0\\
\boldsymbol{\mu}^0
\end{bmatrix}\right\|^2.
\end{equation*}
By inspecting the proof to \cite[Theorem 4.4]{BecTab_FISTA:2009} (and \cite[Lemma
2.3, Lemma 4.1]{BecTab_FISTA:2009}) it is concluded that the optimal
point $\boldsymbol{\lambda}^*,\boldsymbol{\mu}^*$ can be changed to any
feasible point $\boldsymbol{\lambda}_f, \boldsymbol{\mu}_f$ and the
convergence result still holds, i.e.,
\begin{equation*}
D_{N}^{\delta}(\bar{x},\boldsymbol{\lambda}_f,\boldsymbol{\mu}_f)-D_{N}^{\delta}(\bar{x},\boldsymbol{\lambda}^k,\boldsymbol{\mu}^k)\leq
\frac{2L}{(k+1)^2}\left\|\begin{bmatrix}
\boldsymbol{\lambda}_f\\
\boldsymbol{\mu}_f
\end{bmatrix}-\begin{bmatrix}
\boldsymbol{\lambda}^0\\
\boldsymbol{\mu}^0
\end{bmatrix}\right\|^2.
\end{equation*}
That is, there exists a feasible pair
$(\boldsymbol{\lambda}_f,\boldsymbol{\mu}_f)$ such that with finite $k$ we have
\begin{equation}
\label{eq:unboundIters} D_{N}^{\delta}(\bar{x},\boldsymbol{\lambda}^k,\boldsymbol{\mu}^k)> D_{N}^{\delta}(\bar{x},\boldsymbol{\lambda}_f,\boldsymbol{\mu}_f)-\epsilon\ell^*(\bar{x}).
\end{equation}
This implies
\ifoneColumn
\begin{align*}
\delta\mathbf{d}^T\boldsymbol{\mu}^k &\geq
D_{N}^{\delta}(\bar{x},\boldsymbol{\lambda}^k,\boldsymbol{\mu}^k)-V_N^0(\bar{x})
>D_{N}^{\delta}(\bar{x},\boldsymbol{\lambda}_f,\boldsymbol{\mu}_f)-V_N^0(\bar{x})-\epsilon\ell^*(\bar{x})
\geq \epsilon\ell^*(\bar{x})
\end{align*}
\else
\begin{align*}
\delta\mathbf{d}^T\boldsymbol{\mu}^k &\geq
D_{N}^{\delta}(\bar{x},\boldsymbol{\lambda}^k,\boldsymbol{\mu}^k)-V_N^0(\bar{x})\\
&>D_{N}^{\delta}(\bar{x},\boldsymbol{\lambda}_f,\boldsymbol{\mu}_f)-V_N^0(\bar{x})-\epsilon\ell^*(\bar{x})
\geq \epsilon\ell^*(\bar{x})
\end{align*}
\fi
where Lemma~\ref{lem:tighteningRelationsVal} is used in the first
inequality, \eqref{eq:unboundIters} in the second inequality and
\eqref{eq:dualUnbounded} in the final inequality.
This completes the proof.\hfill$\Box$

\subsection{Proof for Theorem~\ref{thm:wellDefStatic}}\label{app:wellDefStatic}
To prove the assertion we need to show that the do loop will exit
for every $\bar{x}\in{\rm{int}}(\mathbb{X}_N^{0})$.
For every point $\bar{x}\in{\rm{int}}(\mathbb{X}_N^{0})$ there exists
$\bar{\delta}\in(0,1)$ such that
$\frac{\bar{x}}{1-\bar{\delta}}\in{\rm{int}}(\mathbb{X}_N^{0})$. Since
${\rm{int}}(\mathbb{X}_N^{0})\subseteq \mathbb{X}_N^{0}$, we have that
$V_N^0(\frac{\bar{x}}{1-\bar{\delta}})<\infty$ and the
optimal solution $\mathbf{y}(\frac{\bar{x}}{1-\bar{\delta}},0)$
satisfies
$\mathbf{A}\mathbf{y}^*(\frac{\bar{x}}{1-\bar{\delta}},0)=\mathbf{b}\frac{\bar{x}}{1-\bar{\delta}}$
and
$\mathbf{C}\mathbf{y}^*(\frac{\bar{x}}{1-\bar{\delta}},0)\leq\mathbf{d}$.
We create the following vector 
\begin{equation}
\label{eq:slaterVec} \bar{\mathbf{y}}(\bar{x}) :=
(1-\bar{\delta})\mathbf{y}^*(\frac{\bar{x}}{1-\bar{\delta}},0)
\end{equation}
which satisfies
\begin{align}
\label{eq:slaterRel1}\mathbf{A}\bar{\mathbf{y}}(\bar{x})&=\mathbf{A}\mathbf{y}^*(\frac{\bar{x}}{1-\bar{\delta}},0)(1-\bar{\delta})=\mathbf{b}\bar{x}\frac{1-\bar{\delta}}{1-\bar{\delta}}=\mathbf{b}\bar{x}\\
\label{eq:slaterRel2}\mathbf{C}\bar{\mathbf{y}}(\bar{x})&=\mathbf{C}\mathbf{y}^*(\frac{\bar{x}}{1-\bar{\delta}},0)(1-\bar{\delta})\leq\mathbf{d}(1-\bar{\delta}).
\end{align}
Hence, by definition \eqref{eq:feasibleRegion} of
$\mathbb{X}_N^{\delta}$ we conclude that for every
$\bar{x}\in{\rm{int}}(\mathbb{X}_N^{0})$ there exist
$\bar{\delta}\in(0,1)$ such that
$\bar{x}\in\mathbb{X}_N^{\bar{\delta}}$. This implies that for every
$\bar{x}\in{\rm{int}}(\mathbb{X}_N^0)$ we have that
either $\bar{x}\in\mathbb{X}_N^{\delta}$ for the current constraint
tightening $\delta\in(0,1)$ or $\bar{x}\notin\mathbb{X}_N^{\delta}$ but
$\bar{x}\in\mathbb{X}_N^{0}$.
Thus, from Lemma~\ref{lem:finiteIters} and Lemma~\ref{lem:detectInfeas} we
conclude that either the
do loop is terminated or $\delta$ is reduced and $l$ is increased for every
$\bar{x}\in{\rm{int}}(\mathbb{X}_N^{0})$ with finite number of
algorithm iterations $k$.

To guarantee that the do loop will terminate for every
$\bar{x}\in{\rm{int}}(\mathbb{X}_N^{0})$, we need to show that the conditions
in the do loop will hold for small enough $\delta$ and with finite
$k$. That is, we need to show that the following two conditions will hold.
\begin{enumerate}
\item For small enough $\delta$, i.e., large enough $l$, we have that
\begin{equation}
\label{eq:dualBoundProof} \delta(\boldsymbol{\mu}^k)^T\mathbf{d}\leq \epsilon\ell^*(\bar{x})
\end{equation}
where $\delta = 2^{-l}\delta_{\rm{init}}$ holds for every algorithm
iteration $k$.
\item For small enough $\delta$, i.e., large enough $l$, the condition
\begin{equation}
\label{eq:perfTestAdaptProof} D_N^{\delta}(\bar{x},\boldsymbol{\lambda}^k,\boldsymbol{\mu}^k) \geq
P_N(A\bar{x}+Bv_0^k,\mathbf{v}_s^k)+\alpha\ell(\bar{x},v_0^k)
\end{equation}
with $\alpha$ satisfying \eqref{eq:alphaBound} holds with finite $k$ whenever 
\begin{equation}
\label{eq:optCondProof} D_N^{\delta}(\bar{x},\boldsymbol{\lambda}^k,\boldsymbol{\mu}^k) \geq
P_N(\bar{x},\mathbf{v}^k)+\frac{\epsilon}{l+1}\ell(\bar{x},v_0^k)
\end{equation}
holds.
\end{enumerate}


We start by showing argument~1.
From the convergence rate of the algorithm \cite{gisAutomatica} it
follows that there exists
$\underline{D}>-\infty$ such that 
$D_{N}^{\delta}(\bar{x},\boldsymbol{\lambda}^k,\boldsymbol{\mu}^k)\geq \underline{D}$ for every algorithm iteration
$k\geq 0$. 
This is used below where we extend the result from
\cite[Lemma 1]{Nedic:2009} to handle the presence of equality
constraints. For algorithm iteration $k\geq 0$,
$\bar{x}\in{\rm{int}}(\mathbb{X}_N^0)$ and $\delta\leq\bar{\delta}/2$ we have
\ifoneColumn
\begin{align*}
\underline{D}&\leq
D_{N}^{\delta}(\bar{x},\boldsymbol{\lambda}^k,\boldsymbol{\mu}^k) 
=\inf_{\mathbf{y}}
\frac{1}{2}\mathbf{y}^T\mathbf{H}\mathbf{y}+(\boldsymbol{\lambda}^k)^T(\mathbf{A}\mathbf{y}-\mathbf{b}\bar{x})+(\boldsymbol{\mu}^k)^T(\mathbf{C}\mathbf{y}-(1-\delta)\mathbf{d})\\
&\leq 
\frac{1}{2}(\bar{\mathbf{y}}(\bar{x}))^T\mathbf{H}\bar{\mathbf{y}}(\bar{x})+(\boldsymbol{\lambda}^k)^T(\mathbf{A}\bar{\mathbf{y}}(\bar{x})-\mathbf{b}\bar{x})+
(\boldsymbol{\mu}^k)^T(\mathbf{C}\bar{\mathbf{y}}(\bar{x})-(1-\delta)\mathbf{d})\\
&\leq
(1-\bar{\delta})^2V_{N}^{0}(\frac{\bar{x}}{1-\bar{\delta}})+(\boldsymbol{\mu}^k)^T(\mathbf{C}\bar{\mathbf{y}}(\bar{x})-(1-\bar{\delta})\mathbf{d})
+(\boldsymbol{\mu}^k)^T\mathbf{d}(\delta-\bar{\delta})\\
&\leq V_{N}^{0}(\frac{\bar{x}}{1-\bar{\delta}})+(\boldsymbol{\mu}^k)^T\mathbf{d}(\delta-\bar{\delta})\\
&\leq V_{N}^{0}(\frac{\bar{x}}{1-\bar{\delta}})-\frac{1}{2}(\boldsymbol{\mu}^k)^T\mathbf{d}\bar{\delta}
\end{align*}
\else
\begin{align*}
\underline{D}&\leq
D_{N}^{\delta}(\bar{x},\boldsymbol{\lambda}^k,\boldsymbol{\mu}^k) \\
&=\inf_{\mathbf{y}}
\frac{1}{2}\mathbf{y}^T\mathbf{H}\mathbf{y}+(\boldsymbol{\lambda}^k)^T(\mathbf{A}\mathbf{y}-\mathbf{b}\bar{x})+\\
&\qquad\qquad\qquad\qquad\qquad\qquad+(\boldsymbol{\mu}^k)^T(\mathbf{C}\mathbf{y}-(1-\delta)\mathbf{d})\\
&\leq 
\frac{1}{2}(\bar{\mathbf{y}}(\bar{x}))^T\mathbf{H}\bar{\mathbf{y}}(\bar{x})+(\boldsymbol{\lambda}^k)^T(\mathbf{A}\bar{\mathbf{y}}(\bar{x})-\mathbf{b}\bar{x})+\\
&\qquad\qquad\qquad\qquad\qquad\qquad+(\boldsymbol{\mu}^k)^T(\mathbf{C}\bar{\mathbf{y}}(\bar{x})-(1-\delta)\mathbf{d})\\
&\leq
(1-\bar{\delta})^2V_{N}^{0}(\frac{\bar{x}}{1-\bar{\delta}})+(\boldsymbol{\mu}^k)^T(\mathbf{C}\bar{\mathbf{y}}(\bar{x})-(1-\bar{\delta})\mathbf{d})+\\
&\qquad\qquad\qquad\qquad\qquad\qquad\qquad\qquad+(\boldsymbol{\mu}^k)^T\mathbf{d}(\delta-\bar{\delta})\\
&\leq V_{N}^{0}(\frac{\bar{x}}{1-\bar{\delta}})+(\boldsymbol{\mu}^k)^T\mathbf{d}(\delta-\bar{\delta})\\
&\leq V_{N}^{0}(\frac{\bar{x}}{1-\bar{\delta}})-\frac{1}{2}(\boldsymbol{\mu}^k)^T\mathbf{d}\bar{\delta}
\end{align*}
\fi
where the equality is by definition, the second inequality holds since
any vector $\bar{\mathbf{y}}(\bar{x})$ is gives larger value than the
infimum, the third and fourth inequalities are due to \eqref{eq:slaterVec}, 
\eqref{eq:slaterRel1} and \eqref{eq:slaterRel2} and since
$(1-\bar{\delta})\in(0,1)$ and the final inequality holds since
$\delta\leq\bar{\delta}/2$. This implies that 
\begin{equation*}
(\boldsymbol{\mu}^k)^T\mathbf{d}\leq
\frac{2(V_{N}^{0}(\frac{\bar{x}}{1-\bar{\delta}})-\underline{D})}{\bar{\delta}}
\end{equation*}
which is finite. We denote by $l_d$ the smallest $l$ such
that $\bar{\delta} \geq 2^{-l_d}\delta_{\rm{init}}$. Since $\delta =
2^{-l}\delta_{\rm{init}}$ this implies that
\ifoneColumn
\begin{align}
\delta (\boldsymbol{\mu}^k)^T\mathbf{d}&\leq
\delta\frac{2(V_{N}^{0}(\frac{\bar{x}}{1-\bar{\delta}})-\underline{D})}{\bar{\delta}}\leq
2^{-l}\delta_{\rm{init}}\frac{2(V_{N}^{0}(\frac{\bar{x}}{1-\bar{\delta}})-\underline{D})}{2^{-l_d}\delta_{\rm{init}}}
\label{eq:dualVarConv} \leq 2^{-l+l_d+1}(V_{N}^{0}(\frac{\bar{x}}{1-\bar{\delta}})-\underline{D})\to 0
\end{align}
\else
\begin{align}
\nonumber\delta (\boldsymbol{\mu}^k)^T\mathbf{d}&\leq
\delta\frac{2(V_{N}^{0}(\frac{\bar{x}}{1-\bar{\delta}})-\underline{D})}{\bar{\delta}}\leq
2^{-l}\delta_{\rm{init}}\frac{2(V_{N}^{0}(\frac{\bar{x}}{1-\bar{\delta}})-\underline{D})}{2^{-l_d}\delta_{\rm{init}}}\\
&\label{eq:dualVarConv} \leq 2^{-l+l_d+1}(V_{N}^{0}(\frac{\bar{x}}{1-\bar{\delta}})-\underline{D})\to 0
\end{align}
\fi
as $l\to\infty$. Especially, with finite $l$ we have that
\eqref{eq:dualBoundProof} holds for every algorithm iteration $k$.
This proves argument~1.

Next we prove argument~2. We start by showing for large enough but finite $l$ that
$P^N(A\bar{x}+B\nu_N(\bar{x}),\mathbf{v}_s^k)$ is finite whenever \eqref{eq:optCondProof}
holds. From the definition of $P_N$ and $\mathbf{v}_s^k$ we have that
$P^N(A\bar{x}+B\nu_N(\bar{x}),\mathbf{v}_s^k)$ is
finite whenever $P_N(\bar{x},\mathbf{v}_s^k)$ is finite and if
$A\xi_{N-1}^k(\bar{x},\delta)\in\mathcal{X}$. For algorithm iteration $k$
such that \eqref{eq:optCondProof} holds we have
\ifoneColumn
\begin{align}
\nonumber \|A(\xi^k_{N-1}(\bar{x},\delta)-z^*_{N-1}(\bar{x},0))\|^2
&\leq
\frac{\|A\|^2}{\lambda_{\min}(H)}\|\xi^k_{N-1}(\bar{x},\delta)-z^*_{N-1}(\bar{x},0)\|_{H}^2\\
\nonumber&\leq
\frac{2\|A\|^2}{\lambda_{\min}(H)}(\delta
(\boldsymbol{\mu}^k)^T\mathbf{d}+\frac{\epsilon}{l+1}\ell^*(\bar{x}))\\
&\leq
\frac{2\|A\|^2}{\lambda_{\min}(H)}\bigg(2^{-l+l_d+1}(V_N^{0}(\frac{\bar{x}}{1-\bar{\delta}})-\underline{D})+\frac{\epsilon}{l+1}\ell^*(\bar{x})\bigg)\to 0
\label{eq:AxiConv}
\end{align}
\else
\begin{align}
\nonumber &\|A(\xi^k_{N-1}(\bar{x},\delta)-z^*_{N-1}(\bar{x},0))\|^2
\leq\\
\nonumber&\leq
\frac{\|A\|^2}{\lambda_{\min}(H)}\|\xi^k_{N-1}(\bar{x},\delta)-z^*_{N-1}(\bar{x},0)\|_{H}^2\\
\nonumber&\leq
\frac{2\|A\|^2}{\lambda_{\min}(H)}(\delta
(\boldsymbol{\mu}^k)^T\mathbf{d}+\frac{\epsilon}{l+1}\ell^*(\bar{x}))\\
&\leq
\frac{2\|A\|^2}{\lambda_{\min}(H)}\bigg(2^{-l+l_d+1}(V_N^{0}(\frac{\bar{x}}{1-\bar{\delta}})-\underline{D})+\frac{\epsilon}{l+1}\ell^*(\bar{x})\bigg)\to 0
\label{eq:AxiConv}
\end{align}
\fi
as $l\to\infty$ where $H={\rm{blkdiag}}(Q,R)$ and the smallest
eigenvalue $\lambda_{\min}(H)>0$ since $H$ is positive
definite. The first inequality follows from Cauchy-Schwarz inequality and
Courant-Fischer-Weyl min-max principle,
the second inequality comes from
Lemma~\ref{lem:tighteningRelationsVar} and the third comes from
\eqref{eq:dualVarConv}. By definition of
$\mathbb{X}_N^{\delta}$ we have $Az_{N-1}^*(\bar{x},0)\in{\rm{int}}(\mathcal{X})$ which
through \eqref{eq:AxiConv} implies that $A\xi_{N-1}^k(\bar{x},\delta)\in\mathcal{X}$
for some large enough by finite $l$, i.e., small enough $\delta$, and
for algorithm iteration $k$ such that \eqref{eq:optCondProof} holds. 

What is left to show is that \eqref{eq:perfTestAdaptProof} 
holds for every
$\alpha\leq 1-2\epsilon-\kappa(\sqrt{2\epsilon}+\sqrt{\Phi_N})^2(\sqrt{2\epsilon}+1)^2$
for large enough but finite $l$ whenever \eqref{eq:optCondProof}
holds. From Lemma~\ref{lem:tighteningRelationsVar} and
\eqref{eq:dualVarConv} we know for 
large enough $l$ and any algorithm iteration $k$ such that
\eqref{eq:optCondProof} holds that
\begin{align*}
\frac{1}{2}\left\|\begin{bmatrix}
\xi^k_{\tau}\\
v^k_{\tau}
\end{bmatrix}-\begin{bmatrix}
z^*_{\tau}\\
v^*_{\tau}
\end{bmatrix}\right\|_{H}^2
&\leq\delta(\boldsymbol{\mu}^k)^T\mathbf{d}+\frac{\epsilon}{l+1}\ell^*(\bar{x})\\
&=2^{-l}\delta_{\rm{init}}(\boldsymbol{\mu}^k)^T\mathbf{d}+\frac{\epsilon}{l+1}\ell^*(\bar{x})\leq 2\epsilon\ell^*(\bar{x})
\end{align*}
for any $\tau = 0,\ldots,N-1$,
where $H = {\rm{blkdiag}}(Q,R)$. Taking the square-root and applying the
reversed triangle inequality gives
\begin{align}
\left|\left\|\begin{bmatrix}
\xi^k_{\tau}\\
v^k_{\tau}
\end{bmatrix}\right\|_{H}-\left\|\begin{bmatrix}
z^*_{\tau}\\
v^*_{\tau}
\end{bmatrix}\right\|_{H}\right|&\leq \left\|\begin{bmatrix}
\xi^k_{\tau}\\
v^k_{\tau}
\end{bmatrix}-\begin{bmatrix}
z^*_{\tau}\\
v^*_{\tau}
\end{bmatrix}\right\|_{H}
\leq 2\sqrt{\epsilon\ell^*(\bar{x})}.
\label{eq:variableBound}
\end{align}
This implies that
\ifoneColumn
\begin{align*}
\left\|\begin{bmatrix}
\xi^k_{N-1}\\
v^k_{N-1}
\end{bmatrix}\right\|_{H}
&\leq \left\|\begin{bmatrix}
z^*_{N-1}\\
v^*_{N-1}
\end{bmatrix}\right\|_{H}+2\sqrt{\epsilon\ell^*(\bar{x})}
=\sqrt{2}\sqrt{\ell(z^*_{N-1},v^*_{N-1})}+2\sqrt{\epsilon\ell^*(\bar{x})}\\
&\leq
\sqrt{2\Phi_N}\sqrt{\ell(z^*_{0},v^*_{0})}+2\sqrt{\epsilon\ell^*(\bar{x})}
\leq (\sqrt{2\Phi_N}+2\sqrt{\epsilon})\sqrt{\ell(z^*_{0},v^*_{0})}\\
&= (\sqrt{\Phi_N}+\sqrt{2\epsilon})\left\|\begin{bmatrix}
z^*_{0}\\
v^*_{0}
\end{bmatrix}\right\|_{H}
\leq (\sqrt{\Phi_N}+\sqrt{2\epsilon})\left(\left\|\begin{bmatrix}
\xi^k_{0}\\
v^k_{0}
\end{bmatrix}\right\|_{H}+2\sqrt{\epsilon\ell^*(\bar{x})}\right)\\
&\leq (\sqrt{\Phi_N}+\sqrt{2\epsilon})(1+\sqrt{2\epsilon})\left\|\begin{bmatrix}
\xi^k_{0}\\
v^k_{0}
\end{bmatrix}\right\|_{H}
\end{align*}
\else
\begin{align*}
\left\|\begin{bmatrix}
\xi^k_{N-1}\\
v^k_{N-1}
\end{bmatrix}\right\|_{H}
&\leq \left\|\begin{bmatrix}
z^*_{N-1}\\
v^*_{N-1}
\end{bmatrix}\right\|_{H}+2\sqrt{\epsilon\ell^*(\bar{x})}\\
&=\sqrt{2}\sqrt{\ell(z^*_{N-1},v^*_{N-1})}+2\sqrt{\epsilon\ell^*(\bar{x})}\\
&\leq
\sqrt{2\Phi_N}\sqrt{\ell(z^*_{0},v^*_{0})}+2\sqrt{\epsilon\ell^*(\bar{x})}\\
&\leq (\sqrt{2\Phi_N}+2\sqrt{\epsilon})\sqrt{\ell(z^*_{0},v^*_{0})}\\
&= (\sqrt{\Phi_N}+\sqrt{2\epsilon})\left\|\begin{bmatrix}
z^*_{0}\\
v^*_{0}
\end{bmatrix}\right\|_{H}\\
&\leq (\sqrt{\Phi_N}+\sqrt{2\epsilon})\left(\left\|\begin{bmatrix}
\xi^k_{0}\\
v^k_{0}
\end{bmatrix}\right\|_{H}+2\sqrt{\epsilon\ell^*(\bar{x})}\right)\\
&\leq (\sqrt{\Phi_N}+\sqrt{2\epsilon})(1+\sqrt{2\epsilon})\left\|\begin{bmatrix}
\xi^k_{0}\\
v^k_{0}
\end{bmatrix}\right\|_{H}
\end{align*}
\fi
where we have used \eqref{eq:variableBound}, $z^*_0=\xi^k_0=\bar{x}$, $\|[z^T v^T]^T\|_{H} =
\sqrt{z^TQz+v^TRv} = \sqrt{2\ell(z,v)}$ and
Definition~\ref{def:ctrbKnown}.
Squaring both sides gives through the definition of $\kappa$ that
\begin{align}
\nonumber\frac{1}{\kappa}\ell^*(A\xi^k_{N-1})&\leq\ell^*(\xi^k_{N-1})=\ell(\xi^k_{N-1},v^k_{N-1})\\
\label{eq:lastFirstStageCostRelation} &\leq (\sqrt{\Phi_N}+\sqrt{2\epsilon})^2(1+\sqrt{2\epsilon})^2\ell(\xi^k_{0},v^k_{0}).
\end{align}
We get for large enough $l$ and for $k$ such that
\eqref{eq:optCondProof} holds that
\ifoneColumn
\begin{align*}
 D_{N}^{\delta}(\bar{x},\boldsymbol{\lambda}^k,\boldsymbol{\mu}^k)&\geq
P_N(\bar{x},\mathbf{v}^k)-\frac{\epsilon}{l+1}\ell^*(\bar{x})
\geq P_N(\bar{x},\mathbf{v}^k)-\epsilon\ell^*(\bar{x})\\
 &=P_N(A\bar{x}+Bv_0^k,\mathbf{v}_s^k)+(1-\epsilon)\ell(\xi^k_{0},v^k_{0})-\ell^*(A\xi^k_{N-1})\\
 &\geq P_N(A\bar{x}+Bv_0^k,\mathbf{v}_s^k)+
\left(1-\epsilon-\kappa(\sqrt{\Phi_N}+\sqrt{2\epsilon})^2(1+\sqrt{2\epsilon})^2\right)\ell(\bar{x},v^k_{0})\\
 &\geq P_N(A\bar{x}+Bv_0^k,\mathbf{v}_s^k)+\alpha\ell(\bar{x},v^k_{0})
\end{align*}
\else
\begin{align*}
 D_{N}^{\delta}&(\bar{x},\boldsymbol{\lambda}^k,\boldsymbol{\mu}^k)\geq\\&\geq
P_N(\bar{x},\mathbf{v}^k)-\frac{\epsilon}{l+1}\ell^*(\bar{x})\\
 &\geq P_N(\bar{x},\mathbf{v}^k)-\epsilon\ell^*(\bar{x})\\
 &=P_N(A\bar{x}+Bv_0^k,\mathbf{v}_s^k)+(1-\epsilon)\ell(\xi^k_{0},v^k_{0})-\ell^*(A\xi^k_{N-1})\\
 &\geq P_N(A\bar{x}+Bv_0^k,\mathbf{v}_s^k)+\\
 &\qquad+
\left(1-\epsilon-\kappa(\sqrt{\Phi_N}+\sqrt{2\epsilon})^2(1+\sqrt{2\epsilon})^2\right)\ell(\bar{x},v^k_{0})\\
 &\geq P_N(A\bar{x}+Bv_0^k,\mathbf{v}_s^k)+\alpha\ell(\bar{x},v^k_{0})
\end{align*}
\fi
where the first inequality comes from
\eqref{eq:optCondProof}, the second
since $l\geq 0$, the equality is due to \eqref{eq:primCostRelation},
the third inequality comes from
\eqref{eq:lastFirstStageCostRelation}, and the final inequality comes
from \eqref{eq:alphaBound}. This concludes the
proof for argument~2. Thus, the do loop will terminate with
finite $l$ and $k$. This implies that $\nu_N$ is defined for
every $\bar{x}\in{\rm{int}}(\mathbb{X}_N^0)$, i.e. that
${\rm{dom}}(\nu_N)\supseteq {\rm{int}}(\mathbb{X}_N^0)$.

Finally, to show \eqref{eq:valueFcnDecreaseThm} we have that
\begin{align*}
V_N^0(\bar{x})&\geq
D_{N}^{\delta}(\bar{x},\boldsymbol{\lambda}^k,\boldsymbol{\mu}^k)-\delta\mathbf{d}^T\boldsymbol{\mu}^k\\
& \geq P_N(A\bar{x}+Bv_0^k,\mathbf{v}_s^k)-\epsilon\ell^*(\bar{x})+
\alpha\ell(\bar{x},v^k_{0})\\
& \geq V_N^0(A\bar{x}+Bv_0^k)+
(\alpha-\epsilon)\ell(\bar{x},v^k_{0})
\end{align*}
where the first inequality comes from
Lemma~\ref{lem:tighteningRelationsVal}, the second from
\eqref{eq:dualBoundProof} and \eqref{eq:perfTestAdaptProof} which
obviously hold also for any $\bar{x}\in{\rm{dom}}(\nu_N)$, and the
third holds since $P_N(A\bar{x}+Bv_0^k,\mathbf{v}_s^k)\geq V_N(A\bar{x}+Bv_0^k)$
and by definition of $\ell^*$. This concludes the proof.
\hfill$\Box$

\else
\fi

\end{document}